\documentclass{gtart_a}
\pdfoutput=1

\usepackage{amscd}

\title[Zero dimensional Donaldson--Thomas invariants of threefolds]{Zero dimensional Donaldson--Thomas\\invariants of threefolds}

\author{Jun Li}
\givenname{Jun}
\surname{Li}
\address{Department of Mathematics\\
Stanford University\\\newline
Stanford, CA 94305\\USA }
\email{jli@math.stanford.edu}
\urladdr{}

\volumenumber{10}
\issuenumber{}
\publicationyear{2006}
\papernumber{46}
\lognumber{0737}
\startpage{2117}
\endpage{2171}

\doi{}
\MR{}
\Zbl{}

\arxivreference{math.AG/0604490}


\received{27 April 2006}
\accepted{10 October 2006}
\proposed{Jim Bryan}
\seconded{Lothar Goettsche, Eleny Ionel}
\revised{}
\published{29 November 2006}
\publishedonline{29 November 2006}
\corresponding{}
\editor{CPR}
\version{}

\let\xysavmatrix\xymatrix
\def\xymatrix{\disablesubscriptcorrection\xysavmatrix}
\AtBeginDocument{\let\bar\wbar\let\tilde\wtilde\let\hat\what}

\DeclareMathOperator{\spec}{Spec} 
 
\DeclareMathOperator{\Ext}{Ext}

 \DeclareMathOperator{\image}{Im}

\DeclareMathOperator{\supp}{Supp}

\renewcommand{\rk}{\mathrm{rk}}
\DeclareMathOperator{\Pic}{Pic}

\DeclareMathOperator{\ext}{\cE\!{\it xt}}

\def\cA{{\mathcal A}}

\def\cE{{\mathcal E}}
\def\cF{{\mathcal F}}

\def\cI{{\mathcal I}}

\def\cN{{\mathcal N}}
\def\cO{{\mathcal O}}
\def\cP{{\mathcal P}}
\def\cQ{{\mathcal Q}}
\def\cR{{\mathcal R}}
\def\cS{{\mathcal S}}

\def\cU{{\mathcal U}}
\def\cV{{\mathcal V}}
\def\cW{{\mathcal W}}

\def\cZ{{\mathcal Z}}

\def\bP{\textbf P}

\def\bA{\textbf A}

\def\kk{{\rm k}}

\def\ZZ{{\mathbb Z}}
\def\QQ{{\mathbb Q}}
\def\RR{{\mathbb R}}
\def\CC{{\mathbb C}}

\def\MM{{\mathfrak M}}
\def\PP{{\mathfrak P}}

\makeatletter
\def\cnewtheorem#1[#2]#3{\newtheorem{#1}{#3}[section]
\expandafter\let\csname c@#1\endcsname\c@prop}

\newtheorem{prop}{Proposition}[section]
\cnewtheorem{theo}[prop]{Theorem}
\cnewtheorem{lemm}[prop]{Lemma}
\cnewtheorem{subl}[prop]{Sublemma}
\cnewtheorem{coro}[prop]{Corollary}
\cnewtheorem{clai}[prop]{Claim}

\theoremstyle{definition}
\cnewtheorem{rema}[prop]{Remark}
\cnewtheorem{exam}[prop]{Example}
\cnewtheorem{defi}[prop]{Definition}
\cnewtheorem{conj}[prop]{Conjecture}
\makeatother  

\def\and{\quad{\rm and}\quad}

\def\bl{\bigl(} \def\br{\bigr)}
\def\Bl{\Bigl(} \def\Br{\Bigr)}

\def\dual{^{\vee}}
\def\ddual{^{\vee\vee}}
\def\Del{\Delta}

\def\defeq{\triangleq}

\def\eps{\epsilon}

\def\kk{\textbf k}

\def\lbe{_{\beta}}
\let\lra=\longrightarrow
\def\lsta{_{\ast}}

\def\lbe{_{\beta}}
\def\lalp{_{\alpha}}

\def\llra{\,\mathop{-\!\!\!\lra}\,}
\def\Lam{\Lambda}

\def\mapright#1{\,\smash{\mathop{\lra}\limits^{#1}}\,}
\def\twomapright#1{\,\smash{\mathop{-\!\!\!\lra}\limits^{#1}}\,}

\let\mh\co

\def\Ob{\cO b}

\def\pr{\text{pr}}
\def\pri{^{\prime}}

\def\Po{{\textbf P}^1}

\def\sub{\subset}
\def\sta{^{\ast}}

\def\upmo{^{-1}}

\def\vir{^{{\rm vir}}}

\def\alphaseq{\alpha_1,\cdots,\alpha_k}

\def\awb{{\alpha\wedge\beta}}
\def\At{{\bA^{\!3}}}

\def\cualp{\cU\hhbr{\alpha}} \def\cvalp{\cV\hhbr{\alpha}}

\def\DT{\mathcal{DT}\!}

\def\hx#1{X^{{[#1]}}}
\def\hhx#1{X^{{[\![#1]\!]}}}  \def\hhy#1{Y^{{[\![#1]\!]}}}

\def\hbr#1{^{{[#1]}}}
\def\hhbr#1{^{{[\![#1]\!]}}}
\def\hhalp{\hhbr{\alpha}} \def\hhbe{\hhbr{\beta}}

\def\halp{^{[\alpha]}} 
\def\hhn{^{[\![n]\!]}}

\def\lalpbe{_{\alpha\beta}}
\def\lralpbe{_{(\alpha,\beta)}}
\def\lrbealp{_{(\beta,\alpha)}}
\def\llam{_{\Lambda}}

\def\pc{_{pc}}
\def\lon{_{[n]}}
\def\lalpi{_{\alpha,i}}
\def\Hil#1{I_{{#1}}(0,n)}

\def\red{\text{red}}

\def\ualp{^\alpha}
\def\ube{^\beta}

\let\virt=\vir

\begin{document}
\begin{abstract}
Using a homotopy approach, we prove in this paper a conjecture of
Maulik, Nekrasov, Okounkov and Pandharipande on the dimension zero
Donaldson--Thomas invariants of all smooth complex threefolds.
\end{abstract}

\maketitle
\setcounter{section}{-1}

\section{Introduction}

Ever since the pioneer work of Donaldson and Thomas on Yang--Mills
theory over Calabi--Yau threefolds \cite{DT,Th}, people have been
searching for their roles in the study of Calabi--Yau geometry and
their relations with other branches of mathematics.  The recent
results and conjectures of Maulik, Nekrasov, Okounkov and
Pandharipande \cite{MNOP1,MNOP2} that relate the invariants of the
moduli of ideal sheaves of curves on a Calabi--Yau manifold to its
Gromov--Witten invariants constitute major progress in this
direction. This paper will address one of their conjectures on
invariants associated to Hilbert scheme of points.

To begin with, we let $(X,H)$ be a smooth projective threefold over
complex numbers $\CC$. For any integer $r\geq 0$, a line bundle
$I\in\Pic(X)$ and two classes $c_2\in H^4(X,\ZZ)$ and $c_3\in
H^6(X,\ZZ)$, we form the moduli space
$$\MM_X^H(r,I,c_2,c_3)
$$
of $H$--stable sheaves of $\cO_X$--modules $\cE$ satisfying
$$\rk \cE=r,\quad \det\cE=I,\quad c_2(\cE)=c_2\and c_3(\cE)=c_3.
$$
Back in the seventies, Maruyama \cite{Maru} proved that such moduli
spaces are quasi-projective, and become projective in case the
quaduple $(r,I,c_2,c_3)$ is relatively prime. Later, Mukai
\cite{Mukai} showed that the first order deformations of any sheaf
$\cE$ in these moduli spaces are given by the traceless part of the
extension group $\Ext^1(\cE,\cE)$; the obstructions to deforming $\cE$
lie in the traceless part of $\Ext^2(\cE,\cE)$ (see also Artamkin
\cite{Atam}). In case $X$ is a Calabi--Yau threefold and $\cE$ is
stable, the traceless part
\begin{gather*}
\Ext^3(\cE,\cE)_0=\Ext^0(\cE,\cE)_0\dual=0\\
\tag*{\hbox{while}}
\Ext^2(\cE,\cE)_0=\Ext^1(\cE,\cE)_0\dual.
\end{gather*}
Hence the moduli space $\MM_X^H(r,I,c_2,c_3)$ admits a
perfect-obstruction theory as defined by Li and Tian \cite{Li-Tian}
and Behrend and Fantechi \cite{Behr-Fant};
when it is projective, it carries a virtual dimension zero cycle
$$[\MM_X^H(r,I,c_2,c_3)]\virt\in H_0\bl\MM_X^H(r,I,c_2,c_3),\ZZ\br.
$$
Its degree is the invariant originally defined and studied by
Donaldson and Thomas \cite{DT}. Following \cite{MNOP1,MNOP2}, we shall call
them Donaldson--Thomas invariants of the Calabi--Yau manifold $X$.

One special class of such moduli space studied extensively in
\cite{MNOP1,MNOP2} is when $r=1$ and $I=\cO_X$. Because $X$ is smooth and
$\cE$ is stable, it must be torsion free, and be a subsheaf of its
double dual $\cE\ddual\cong\cO_X$; hence it becomes an ideal sheaf
of a subscheme $Z\sub X$. Following the notation of \cite{MNOP1,MNOP2},
after picking a curve class $\beta$ and an integer $n$, we denote by
$$I_X(\beta,n)$$
the Hilbert scheme of one dimensional subschemes $Z\sub X$
satisfying $[Z]=\beta$ and $\chi(\cO_Z)=n$. A simple argument shows
that $I_X(\beta,n)$ is the moduli space $\MM_X(0,I, \beta, c_3)$
with $c_3$ the third Chern class of any ideal sheaf of $Z\sub X$ in
$I_X(\beta,n)$. One special feature of the moduli of rank one
torsion free sheaves is that they admit perfect-obstruction theory
for all smooth projective threefolds.

Following \cite{MNOP1,MNOP2}, for Calabi--Yau threefold $X$ and curve
$\beta$ one forms the generating function
$$\DT_{X,\beta}(q)=\sum_n \deg[I_X(\beta,n)]\virt q^n.
$$
In case $X$ is any smooth threefold, since the Hilbert schemes
$I_X(0,n)$ have virtual dimensions zero, one defines
$\DT_{X,\beta}(q)$ according to the same formula and call it the
dimensional zero Donaldson--Thomas series as well. Of the several
conjectures on $\DT_{X,\beta}(q)$ proposed in \cite{MNOP1,MNOP2}, one is
about the dimension zero Donaldson--Thomas invariants $\DT_{X,0}(q)$.
Let $M(q)$ be the three dimensional partition function
$$M(q)=\prod_n\frac{1}{(1-q^n)^n}
$$
and $c_3(T_X\otimes K_X)$ be the third Chern class, viewed as the
Chern number, of $T_X\otimes K_X$.

\begin{conj}{\rm\cite{MNOP1,MNOP2}}\label{mnop}\qua 
For any smooth projective threefold $X$, the
dimension zero Donaldson--Thomas series $\DT_{X,0}(q)$ has the form
$$\DT_{X,0}(q)=M(-q)^{c_3(T_X\otimes K_X)}.
$$
\end{conj}

In this paper, we shall prove this conjecture for all compact smooth
complex threefolds.

\begin{theo}\label{theorem}
The zero-dimensional Donaldson--Thomas series $\DT_{X,0}(q)$ for any
compact smooth complex threefold $X$ are of the form
$$\DT_{X,0}(q)=M(-q)^{c_3(T_X\otimes K_X)}.
$$
\end{theo}

We remark that this conjecture was independently proved for the class
of Calabi--Yau threefolds based on different method by Behrend and
Fantechi \cite{CY}, and for projective threefolds by Levine and
Pandharipande \cite{LP}.

We now briefly outline the proof of this theorem. Clearly, in case
$X$ is a disjoint union of two smooth proper threefolds $X_1$ and
$X_2$, then
$$I_X(0,n)=\coprod_{n_1+n_2=n} I_{X_1}(0,n_1)\times
I_{X_2}(0,n_2).
$$
Since the deformation of the ideal sheaf of the union $Z_1\cup
Z_1\sub X_1\cup X_2$ is the direct product of the deformation of
$Z_1\sub X_1$ and the deformation of $Z_2\sub X_2$, we have
$$[I_{X}(0,n)]\virt=\coprod_{n_1+n_2=n} [I_{X_1}(0,n_1)]\virt\times
[I_{X_2}(0,n_2)]\virt.
$$
Therefore
$$\DT_{X,0}(q)=\DT_{X_1,0}(q)\cdot \DT_{X_2,0}(q).
$$
Put differently, the correspondence that sends any threefold to its
zero-dimensional Donaldson--Thomas series defines a homomorphism from
the additive semigroup
$$\PP_\CC=\bigl\{\text{All smooth projective
threefolds}\bigr\}/\text{iso}
$$
to the multiplicative semigroup of infinite series $\ZZ[[q]]$.

A distant cousin of complex manifolds are so called weakly complex
manifolds, which by definition are smooth compact real manifolds $M$
(possibly with smooth boundaries) together with $\CC$--vector bundle
structures on the stabilizations\footnote{Here as usual a
stabilization of $T_M$ is a direct sum of $T_M$ with some trivial
bundle $\RR^m$ on $M$; together with the identity
$T_M\oplus\RR^m=(T_M\oplus \RR^m)\oplus \RR^2$ with the last $\RR^2$
is given the obvious complex structure $\RR^2\cong\CC$.} of their
tangent bundles $T_M$.

The equivalence classes of weakly complex manifolds (without
boundaries) modulo the relations $[\partial W]=0$ forms a group,
called the complex cobordism group $\Omega^\CC$, under the addition
$[M_1]+[M_2]=[M_1\coprod M_2]$. It is a classical result that
$\Omega^\CC\otimes_\ZZ\QQ$ is generated by all possible products of
projective spaces $\bP^n$. As a consequence, the six (real)
dimensional complex cobordism group $\Omega^\CC_6$ is generated by
$$Y_1=\bP^3,\quad Y_2=\bP^2\!\times\!\Po \and Y_3=(\Po)^3.
$$
The crucial step in proving \fullref{theorem} is to establish

\begin{prop}\label{prop}
There are universal polynomials
$$f_0,f_1,f_2,\cdots
$$
in Chern numbers of smooth complex threefolds so that for any smooth
projective threefold $X$ its zero-dimensional Donaldson--Thomas
series are of the form
$$\DT_{X,0}(q)=\sum_n f_n(X)q^n.
$$
\end{prop}

Because any complex threefold $X$ is $\CC$--cobordant to
$\frac{m_1}{m}Y_1+\frac{m_2}{m}Y_2+\frac{m_3}{m}Y_3$ for some
integers $n_i$, knowing the Proposition, and that cobordant weakly
complex manifolds have identical Chern numbers,
$$\sum_n f_n(mX)q^n=\sum_n f_n(m_1Y_1+m_2Y_2+m_3Y_3)q^n.
$$
By the definition, the left hand side is
$$\DT_{mX,0}(q)=\DT_{X,0}(q)^m
$$
while the right hand side is
$$\DT_{m_1Y_1+m_2Y_2+m_3Y_3,0}(q)=\DT_{Y_1,0}(q)^{m_1}\cdot
\DT_{Y_2,0}(q)^{m_2}\cdot \DT_{Y_3,0}(q)^{m_3}.
$$
Because $Y_i$ are toric threefolds, their Donaldson--Thomas series
are known \cite{MNOP1,MNOP2} to have the form
$$\DT_{Y_i,0}(q)=M(-q)^{c_3(T_{Y_i}\otimes K_{Y_i})}.
$$
Put together, and adding that
$$m\, c_3(T_{X}\otimes K_{X})=\sum_{i=1}^3 m_i\, c_3(T_{Y_i}\otimes
K_{Y_i}),
$$
we obtain
$$\DT_{X,0}(q)^m=M(-q)^{\sum_{i=1}^3 m_i\, c_3(T_{Y_i}\otimes
K_{Y_i})}=M(-q)^{ m\, c_3(T_{X}\otimes K_{X})}.
$$
Finally, because both $\DT_{X,0}(q)$ and $M(-q)$ are power series
with integer coefficients and constant coefficient one,
$$\DT_{X,0}(q)=M(-q)^{c_3(T_X\otimes
K_X)}.
$$
This would prove the theorem.

As to the proof of the Proposition, we shall first construct a
collection of approximations $\hx{\tau}$ of $\hx{n}$ indexed by
partitions of $[n]=\{1,\cdots,n\}$. To each $\hx{\tau}$, we shall
construct its virtual cycle and prove that its degree can be
approximated by the degree of the virtual cycle of $\hx{\sigma}$ of
$\sigma<\tau$, with errors expressible in terms of a universal
expression of the Chern numbers of $X$. Thus by induction, we prove
that the degree of the virtual cycle of $\hx{n}$ can also be
expressed universally in terms of the Chern numbers of $X$, thus
proving the Proposition and the Theorem.

During the early stage of this work, K~Behrend developed a theory
of micro-local analysis for a symmetric obstruction theory
\cite{symm}; later, jointly with B~Fantechi they proved \fullref{mnop} for Calabi--Yau threefolds \cite{CY}. Toward the end of
finalizing this paper, the author was kindly informed by
R~Pandharipande that he and M~Levine have proved the same
conjecture using algebraic K--theory \cite{LP}.

The author also like to take this opportunity to thank Weiping Li
for his valuable comments.

\subsection{Terminology}

We shall work with the category of analytic functions in this paper.
Thus for any reduced quasi-projective scheme $W$, we shall work with
and denote by $\cO_W$ the sheaf of analytic functions on $W$; we
shall use ordinary open subsets of the $W$ unless otherwise stated.

To distinguish from that of analytic spaces, we shall reserve the
words schemes and morphisms to mean algebraic schemes and algebraic
morphisms, though viewed as objects in analytic category.

Often, we shall work with open subsets of quasi-projective schemes.
We shall endow such sets with their reduced induced analytic
structures and work with their sheaves of analytic functions. We
shall call such space either analytic spaces or analytic schemes.
Accordingly, whenever we say two analytic spaces isomorphic we mean
that they are isomorphic as analytic spaces.

In this paper, we reserve the word smooth to mean the smoothness in
$C^\infty $--category. Thus a smooth function is a
$C^\infty$--function and a smooth map is a $C^\infty$--map.

There is one case in which we need to keep non-reduced scheme
structures. It is the case of flat $U$--families of zero-subschemes
$\cZ\sub U\times Y$ with $U$ an analytic space and $Y$ an open
subset of a smooth varieties. In this case, $\cZ$ is defined by the
ideal sheaf $\cI_\cZ\sub\cO_{U\times Y}$; the structure sheaf of
$\cZ$ is the quotient $\cO_{U\times Y}/\cI_\cZ$; we say that $\cZ$
is flat over $U$ if $\cO_{\cZ}$ is flat over $\cO_U$.

\section{Hilbert schemes of $\alpha$--points}

The purpose of this section is to construct the filtered
approximation of $I_X(0,n)$ by Hilbert schemes of $\alpha$--points.

For convenience, we let
$$X\hbr{n}=I_X(0,n)_\red
$$
be the Hilbert scheme $I_X(0,n)$ with the reduced scheme structure.

\subsection{The Definition}

We begin with a finite set $\Lambda$ of order $|\Lambda|$; it could
be the set of $n$ integers $[n]=\{1,\cdots,n\}$ or a subset of
$[n]$. For such $\Lam$ we follow the convention
$$X^\Lambda=\{(x_a)_{a\in\Lam}\mid x_a\in X\}.
$$
In case $|\Lam|=n$, we define
$$X^{(\Lam)}=X^{(n)}=S^n X\and
X^{[\Lam]}=X^{[n]}.
$$
Using the Hilbert--Chow morphism $hc\mh\hx{\Lambda}\to X^{(\Lam)}$,
we define
$$\hhx{\Lambda}=\hx{\Lambda}\times_{X^{(\Lambda)}}X^\Lambda;
$$
it comes with a tautological projection $\hhx{\Lambda}\to
X^\Lambda$. Obviously, according to this definition a closed point
of $X\hhbr{\Lam}$ is a pair of a $0$--scheme $\xi\in X^{[\Lam]}$ and
a point $(x_a)_{a\in\Lam}\in X^\Lam$ such that $hc(\xi)=\sum x_a$.

We next look at the set $\cP_\Lambda$ of all partitions of, or
equivalence relations on, the set $\Lambda$. In case
$\alpha\in\cP_\Lambda$ has $k$ equivalence classes
$\alpha_1,\cdots,\alpha_k$, we write
\begin{equation}\label{alpha}
\alpha=(\alpha_1,\cdots,\alpha_k).
\end{equation}
It has two distinguished elements: one is $\Lam$ that is the
partition with a single equivalence class $\Lam$; the other is
$0_\Lam$ that is the partition whose equivalence classes are all
single element sets. The set $\cP_\Lambda$ has a partial ordering
``$\leq$'' defined by
$$``\alpha\geq\beta"\Longleftrightarrow ``a\sim_\beta b \Rightarrow
a\sim_\alpha b".
$$
In case $\alpha=(\alpha_1,\cdots,\alpha_k)$, then $\alpha\geq\beta$
if and only if $\beta$ is finer than $\alpha$, or that $\beta$ can
be written as
$$\beta=(\beta_{11},\cdots,\beta_{1l_1},\cdots,
\beta_{k1},\cdots,\beta_{kl_k})\quad\text{so that}\quad
\alpha_i=\cup_{j=1}^{l_i}\beta_{ij}.
$$
Under this partial ordering, the element $\alpha\wedge\beta$, which
is defined by
$$``a\sim_{\alpha\wedge\beta}b"\Longleftrightarrow ``a\sim_\alpha b
\ {\rm and}\ a\sim_\beta b",
$$
or equivalently for $\beta=(\beta_1,\cdots,\beta_l)$ it is
$\alpha\wedge\beta=(\alpha_1\cap\beta_1,\cdots,\alpha_k\cap\beta_l)$,
is the largest element among all that are less than or equal to both
$\alpha$ and $\beta$. Following this rule, $0_\Lam$ is the smallest
element and $\Lam$ is the largest element in $\cP_\Lam$.

For $\alpha\in\cP_\Lambda$, we let $\hhx{\alpha}$ be (reduced)
Hilbert scheme of $\alpha$--points in $X$. Let $\alpha$ be as in
\eqref{alpha}. Because each $\alpha_i$ is a set, we can form
$X^{\alpha_i}$, $X^{(\alpha_i)}$, $X^{[\alpha_i]}$ and
$X\hhbr{\alpha_i}$ respectively. We then define
$$
X^{(\alpha)}=\prod_{i=1}^k X^{(\alpha_i)},\quad
\hx{\alpha}=\prod_{i=1}^k \hx{\alpha_i},\quad
\hhx{\alpha}=\prod_{i=1}^k \hhx{\alpha_i};
$$
they fit into the Cartesian product:
$$\begin{CD}
\hhx{\alpha} @>>> \hx{\alpha}\\
@VVV @VV{hc}V\\
X^\Lam @>{S\lalp}>> X^{(\alpha)}
\end{CD}
$$
We shall call points in $X\hhalp$ $\alpha$--zero-subschemes and call
$X\hhalp$ the Hilbert scheme of $\alpha$--points.

The space $X\hhalp$ coincides with $X\hhbe$ over a large open subset
of each. For instance, both $X\hhbr{0_\Lam}=X^\Lam$ and
$X\hhbr{\Lam}$ contains as their open subsets the set of $n=|\Lam|$
distinct ordered points in $X$. It is when distinct simple points
specialize to points with multiplicities the space $X\hhbr{0_\Lam}$
becomes different from $X\hhbr{\Lam}$: for the former they remain as
simple points by allowing multiple simple points to occupy identical
positions in $X$; for the later fat points with non-reduced scheme
structures emerge. This way, the collection
$$\{\hhx{\beta}\mid\beta\in\cP_\Lambda\}
$$
forms an increasingly finer approximation of $\hhx{\Lambda}$.

Though the notion of $X\hhbr{\alpha}$ seems artificial at first, it
proves to be useful in keeping tracking of the difference among all
$X\hhbr{\alpha}$. Lastly, in case $\Lam=[n]$, we shall follow the
convention
$$X^\Lam=X^n,\quad X^{(\Lam)}=X^{(n)},\quad
X\hbr{\Lam}=X\hbr{n},\quad X\hhbr{\Lam}=X\hhbr{n}.
$$

\subsection{The relative case}

We can generalize the notion of Hilbert scheme of $\alpha$--points to
that of smooth families of varieties. Let
$$\pi\co  Y\lra T
$$
be a smooth family of quasi-projective varieties. For any integer
$l$, we let $I_{Y/T}(0,l)$ be the relative Hilbert scheme of length
$l$ $0$--subschemes of fibers of $Y/T$. It is the fine moduli scheme
representing the functor parameterizing all flat $S$--families of
length $l$ $0$--schemes $Z\sub Y\times_T S$. As before, we let
$$Y\hbr{l}=I_{Y/T}(0,l)_\red
$$
be $I_{Y/T}(0,l)$ with the reduced scheme structure. The moduli
$Y\hbr{l}$ is a scheme over $T$ with a universal family. In case for
an open $U\sub T$ with $Y\times_T U=Y_0\times U$, then canonically
$$Y\hbr{l}\times_TU\equiv  Y_0\hbr{l}\times U.
$$
As before, for any $\alpha\in\cP_\Lam$ we let $Y^\Lam$,
$Y^{(\alpha)}$, $Y^{[\alpha]}$ and $Y\hhbr{\alpha}$ be the Cartesian
products presented before with $X$ replaced by $Y$ and with products
replaced by fiber products over the base scheme $T$.
Again in case $Y\times_T U=Y_0\times U$,
$$Y\hhbr{\alpha}\times_T U\equiv  Y_0\hhbr{\alpha}\times U.
$$
The three cases we shall apply this construction is for the trivial
fiber bundle $\pr_1\mh X\times X\to X$, for the total space of the
tangent bundle $TX\to X$ and for the total space of the universal
quotient bundle $Q$ of a Grassmannian $Gr=Gr(N,3)$ of quotients
$\CC^3$ of $\CC^N$. We shall come back to this in detail later.

\subsection{Partial equivalences}\label{sec1.3}

As mentioned before, the collection $\hhy{\beta}$ forms an
increasingly finer approximation of $\hhy{\Lambda}$. It is the
purpose of this subsection to make this precise.

We begin with comparing $Y\hhbr{\Lam}$ with $Y\hhbr{\alpha}$ for an
$\alpha=(\alpha_1,\cdots,\alpha_k)$. By definition, a point in
$Y\hhbr{\alpha_i}$ consists of
$$(\xi_i,(x_a)_{a\in\alpha_i})\in Y\hbr{\alpha_i}\times_T Y^{\alpha_i}
$$
subject to the constraint $hc(\xi_i)=\sum_{a\in\alpha_i} x_a$. In
case the support of $\xi_i$ is disjoint from that of $\xi_j$, then
$\xi_i\cup\xi_j$ is naturally a zero-subscheme in $Y$ of length
$|\alpha_i\cup\alpha_j|$; the pair
$(\xi_i\cup\xi_j,(x_a)_{a\in\alpha_i\cup\alpha_j})$ thus is a point
in $Y\hhbr{\alpha_i\cup\alpha_j}$. Applying this to all pairs $1\leq
i<j\leq k$, we see that
$$\bl\cup_{i=1}^k\xi_i, (x_a)_{a\in\Lam}\br\in Y\hhbr{\Lam}
$$
if and only if the supports $hc(\xi_i)$ are mutually disjoint.

In general, for $\alpha<\beta=(\beta_1,\cdots,\beta_l)$ and $1\leq
j\leq l$, the pair
$$\bl \cup_{\alpha_i\sub\beta_j}\xi_i,(x_a)_{a\in\beta_j}\br\in
Y\hhbr{\beta_j}
$$
if and only if the supports $\{hc(\xi_i)\mid \alpha_i\sub\beta_j\}$
are mutually disjoint. This leads to the definition

\begin{defi}
For $\alpha<\beta$ we define $\Delta_{(\alpha,\beta)}$ be the set
$$\{x\in Y^\Lambda\mid x_a=x_b\ \text{for at least one pair $a,b\in
\Lambda$ so that $a\sim_\beta b$ but $a\not\sim_{\alpha} b$}\,\};
$$
for general $\alpha\ne\beta$ we define
$\Delta_{(\alpha,\beta)}=\Delta_{(\alpha,\alpha\wedge\beta)}
\cup\Delta_{(\beta,\alpha\wedge\beta)}$; we define the discrepancy
between $Y\hhbr{\alpha}$ and $Y\hhbr{\beta}$ be
$$\Delta_{(\alpha,\beta)}\hhbr{\alpha}\defeq
\hhy{\alpha}\times_{Y^\Lambda}\Delta_{(\alpha,\beta)}.
$$
\end{defi}

$$\hhy{\alpha}_{(\alpha,\beta)}=\hhy{\alpha}-\Delta\hhbr{\alpha}_{(\alpha,\beta)}.
\leqno{\hbox{We define}}$$

\begin{lemm}\label{Hil-iso}
Given any pair $\alpha,\beta\in\cP_\Lambda$, we have a functorial
isomorphism
$$\hhy{\alpha}_{(\alpha,\beta)}\cong
\hhy{\beta}\lrbealp.
$$
\end{lemm}

Here by functorial isomorphisms we mean those that are induced by
the universal property of the respective moduli spaces.

\begin{proof}[Proof of \fullref{Hil-iso}]
We first prove the case $\alpha=\Lambda$. Let
$\beta=(\beta_1,\cdots,\beta_l)\in\cP_\Lambda$ with $m_i=|\beta_i|$;
let $S=\hhy{\Lambda}_{(\Lambda,\beta)}$; and let $(\cW,\varphi)$ be
the tautological family of $S$. By definition, $\varphi\mh S\to
Y^\Lam$ is the tautological map and $\cW$ is a flat $S$--family of
length $|\Lam|$ zero-subschemes in $Y$ whose Hilbert--Chow map
$hc_{\cW}\mh S\to Y^{(\Lam)}$ coincide with the composite of
$\varphi$ and $Y^\Lam\to Y^{(\Lam)}$.

Similarly, a $T$--morphism $S\to\hhy{\beta}$ 
is classified by a collection
\begin{equation}
\label{S}\bl\cW_1,\cdots,\cW_l;\varphi_{1},\cdots,\varphi_{l}\br
\end{equation}
of $S$--families $\cW_i\sub Y\times_T S$ in $\hx{\beta_i}$ and
$\varphi_{i}\mh S\to Y^{\beta_i}$ so that the Hilbert--Chow morphism
$hc_{\cW_i}\mh S\to \hx{\beta_i}$ (of the family $\cW_i$) makes the
diagram
\begin{equation}\label{sq1}
\begin{CD}
S @>{(\rho_{\cW_1},\cdots,\rho_{\cW_l})}>> \hx{\beta}\\
@V{(\varphi_1,\cdots,\varphi_l)}VV @VVV\\
Y^\Lam @>{\quad\qquad\qquad}>> Y^{(\beta)}\\
\end{CD}
\end{equation}
commutative.

To proceed, we shall show that we can split the family
$(\cW,\varphi)$ into a family as in \eqref{S}. For each $a\in
\Lambda$, we denote by $\varphi_a\mh S\to Y_a$ the $a$-th component
of $\varphi$, by $\Gamma_{\varphi_a}\sub Y\times_T S$ the graph of
$\varphi_a$ and by $\Gamma_{\beta_i}$ the union $\cup_{a\in
\beta_i}\Gamma_{\varphi_a}$. By the definition of
$\hhy{\Lambda}_{(\Lambda,\beta)}$, the set
$\Gamma_{\beta_1},\cdots,\Gamma_{\beta_l}$ are mutually disjoint
closed subsets of $Y\times S$; hence are open and closed subsets of
$\Gamma=\cup_{a\in \Lambda} \Gamma_{\varphi_a}$.

On the other hand, the closed subscheme $\iota\mh \cW\hookrightarrow
Y\times_T S$ is set-theoretically identical to $\Gamma$; hence
$\cW_i=\iota\upmo(\Gamma_{\beta_i})$ are mutually disjoint open and
closed subsets of $\cW$, which therefore inherit scheme structures
from $\cW$ so that $\cW=\coprod_{i=1}^l\cW_i$. In particular, they
are close subschemes of $Y\times_T S$, flat and finite over $S$.
Finally, because the fibers of $\cW_i$ over general closed points in
$S$ have length $m_i$, by the flatness, $\cW_i$ is a family of
length $m_i$ zero-subschemes in $Y$.

Now, because the sets $\Gamma_{\beta_i}$ are mutually disjoint, the
associated classifying morphisms
$$\rho_{\cW_i}\co  S\lra Y\hbr{\beta_i} \and
\varphi_{i}=\prod_{a\in \beta_i}\varphi_a\co  S\lra Y^{\beta_i}
$$
satisfies the commutative diagram \eqref{sq1}. Therefore, each pair
$(\cW_i, \varphi_{i})$ associates to a canonical morphism
$$\eta_{i}\co S\lra Y\hbr{\beta_i}\times_{Y^{(\beta_i)}} Y^{\beta_i}
=\hhy{\beta_i}.
$$
Put them together, we obtain a morphism
$$\eta=(\eta_1,\cdots,\eta_l)\co S\lra
\hhy{\beta_1}\times_T\cdots\times_T\hhy{\beta_l}=\hhy{\beta}
$$
whose image is contained in $\hhy{\beta}_{(\beta, \Lambda)}$. This
way, we have constructed an induced morphism
$$\eta\co  S=\hhy{\Lambda}_{(\Lambda,\beta)}\lra \hhy{\beta}_{(\beta, \Lambda)}.
$$
To complete the proof of this special case, we need to construct a
morphism
$$\eta\co  \hhy{\beta}_{(\beta, \Lambda)} \lra S=
\hhy{\Lambda}_{(\Lambda, \beta)}
$$
that is the inverse of $\eta$. But this is straight forward and
shall be omitted. This proves the lemma for the case
$\alpha=\Lambda$.

For the general case $\alpha=(\alphaseq)\in\cP_\Lambda$, the
$\alpha\wedge\beta$ consists of equivalence classes $\alpha_i\cap
\beta_j$, each of order $m_{ij}$. Obviously,
$$\alpha_i\wedge\beta\triangleq (\alpha_i\wedge\beta_{1},\cdots,\alpha_i\wedge\beta_{l})
\in\cP_{\alpha_i}.
$$
Hence we can apply the proven case of this lemma to conclude
$$\hhy{\alpha_i\wedge\beta}_{(\alpha_i\wedge\beta,\alpha_i)}\cong
\hhy{\alpha_i}_{(\alpha_i,\alpha_i\wedge\beta)}.
$$
Therefore, because $\hhy{\beta}=\prod_T^k \hhy{\beta_i}$, one checks
easily that:
\begin{align*}
\hhy{\alpha}_{(\alpha, \alpha\wedge\beta)}
&=\prod_{i=1}^k \hhy{\alpha_i}_{(\alpha_i,\alpha_i\wedge\beta)} \cong
\prod_{i=1}^k
\hhy{\alpha_i\wedge\beta}_{(\alpha_i\wedge\beta,\alpha_i)}\\
&=\prod_{i=1}^k\Bigl(\prod_{j=1}^l
\hhy{\alpha_i\wedge\beta_{j}}\Bigr)\times_{Y^{\alpha_i}} \bl
Y^{\alpha_i}-\Delta_{(\alpha_i,\alpha_i\wedge\beta)}\br\\
&=\Bigl(\prod_{i,j}
\hhy{\alpha_i\wedge\beta_{j}}\Bigr)\times_{Y^\Lambda} \bl
Y^{\Lambda}-\Delta_{(\alpha,\awb)}\br
\end{align*}
For the same reason,
$$\hhy{\beta}_{(\beta,\alpha\wedge\beta)}\cong
\Bigl(\prod_{i,j}
\hhy{\alpha_i\wedge\beta_{j}}\Bigr)\times_{Y^\Lambda} \bl
Y^\Lambda-\Delta_{(\beta,\alpha\wedge\beta)}\br.
$$
Because
$\Delta_{(\alpha,\beta)}=\Delta_{(\alpha,\alpha\wedge\beta)}\cup
\Delta_{(\beta,\alpha\wedge\beta)}$, we obtain
$$\hhy{\alpha}_{(\alpha,\beta)}=\hhy{\alpha}_{(\alpha,\alpha\wedge\beta)}
\times_{Y^\Lambda}\bl Y^\Lambda-\Delta_{(\alpha,\beta)}\br\cong
\hhy{\beta}_{(\beta,\alpha\wedge\beta)} \times_{Y^\Lambda}\bl
Y^\Lambda-\Delta_{(\alpha,\beta)}\br=\hhy{\beta}_{(\alpha,\beta)}.
$$
This proves the Lemma.
\end{proof}

\subsection{Universal families under partial equivalence}

Let $(\alpha,\beta)$ and
$\alpha\wedge\beta=(\alpha_1\cap\beta_1,\cdots,\alpha_k\cap\beta_l)$
be as in the proof of the previous lemma; let
$$(\cW_{ij},\varphi_{ij}: i=1,\cdots,k, j=1,\cdots,l)
$$
be the universal family of $\hhy{\alpha\wedge\beta}$. Under the
partial equivalence $
\hhy{\alpha}_{(\alpha,\alpha\wedge\beta)} \cong
\hhy{\alpha\wedge\beta}_{(\alpha\wedge\beta,\alpha)}$, the
restriction of $\cW_{ij}$ to $\hhy{\awb}_{(\awb,\alpha)}$:
$$\cW_{i1}\times_{\hhy{\awb}}\hhy{\awb}_{(\awb,\alpha)},\cdots,
\cW_{il}\times_{\hhy{\awb}}\hhy{\awb}_{(\awb,\alpha)},
$$
are families of zero-schemes of $Y$ over
$\hhy{\alpha}_{(\alpha\wedge\beta,\alpha)}$. Following the previous
proof, these families, viewed as subschemes in
$Y\times_T\hhy{\awb}_{(\awb,\alpha)}$, are mutually disjoint. Hence
their union
$$\Bigl(\coprod_{j=1}^l \cW_{ij}\Bigr)
\times_{\hhy{\awb}}\hhy{\awb}_{(\awb,\alpha)}
$$
forms a flat family of zero schemes in $Y$ of length $|\alpha_i|$
over $\hhy{\awb}_{(\awb,\alpha)}$.

\begin{coro}\label{cor1.2}
Let $(\cZ_i,\varphi_i; i=1,\cdots,k)$ be the universal family of
$Y\hhbr{\alpha}$. Then for each $i$,
$$\Bigl(\coprod_{j=1}^l \cW_{ij}\Bigr)
\times_{\hhy{\awb}}\hhy{\awb}_{(\awb,\alpha)} =\!\!\!= \cZ_i
\times_{\hhy{\alpha}}\hhy{\alpha}_{(\alpha,\awb)}
$$
as families of relative zero-subschemes in $Y/T$.
\end{coro}

\subsection{Hilbert scheme of centered $\alpha$-points}

In order to parameterize family of slices in $Y\hhbr{\alpha}$, we
need the notion of Hilbert scheme of centered $\alpha$--points.

Let $\pi\mh Y\to T$ be the total space of a rank three vector
bundle, viewed as a smooth family of affine schemes isomorphic to
$\At$. We let
$$\varpi\co  Y^\Lam\lra V, \ (x_a)_{a\in\Lam}\in Y_t^\Lam\longmapsto
\frac{1}{|\Lam|}\sum_{a} x_a\in Y_t
$$
be the fiberwise averaging morphism and let
$$\varpi\hhbr{\alpha}\co  Y\hhbr{\alpha}\llra Y
$$
be its composition with the tautological $Y\hhalp\to Y^\Lam$. We
define the relative Hilbert scheme of centered $\alpha$--points be
the preimage of the zero section $0_Y$ of $Y$ under
$\varpi\hhbr{\alpha}$:
\begin{equation}
\label{20.1} Y_0\hhbr{\alpha}=Y\hhbr{\alpha}\times_Y 0_Y.
\end{equation}
Intuitively, $Y\hhalp_0$ consists of $\alpha$--zero-subschemes whose
center of support lie in the zero section of $Y$.

For any pair $\alpha,\beta\in\cP_\Lam$, we define
$$Y_{0,(\alpha,\beta)}\hhbr{\alpha}=Y_{0}\hhbr{\alpha}\cap
Y_{(\alpha,\beta)}\hhbr{\alpha}.
$$
The partial equivalence for $Y\hhalp$ carries over to
\begin{equation}
\label{par2}
Y_{0,(\alpha,\beta)}\hhbr{\alpha}=Y_{0,(\beta,\alpha)}\hhbr{\beta}
\end{equation}

\subsection{Outline of the proof}

We now explain briefly the strategy to prove the main Proposition.
First, because $X\hhbr{n}$ is finite over $X\hbr{n}$, the degree of
its virtual cycle is a fraction of that of $X\hbr{n}$. To study the
former, we first make sense of the virtual cycles
$[X\hhbr{\alpha}]\virt$ for all partitions $\alpha$ of $[n]$; we
then construct explicitly their cycle representatives $D\lalp$ as
cycles in $X^n$. Using the isomorphism
$X\hhbr{\alpha}_{(\alpha,\beta)}\cong
X\hhbr{\beta}_{(\beta,\alpha)}$, we can choose $D\lalp$ and
$D_\beta$ so that their difference lies entirely in a small tubular
neighborhood of $\Delta_{(\alpha,\beta)}\sub X^n$. Repeating this
procedure, we show that the desired degree $\deg D_{[n]}$ is a
linear combination of $\deg D\lalp$ plus a discrepancy term which we
denote by $\delta_{[n]}$. Using induction on $n$, to prove the main
proposition we only need to show that $\delta_{[n]}$ only depend on
the Chern numbers of $X$.

To prove the last statement, we shall find a cycle representative of
$\delta_{[n]}$ that is entirely contained in a small (tubular)
neighborhood of the top diagonal
$$X^\Lam_\Del=\{(x,\cdots,x)\mid x\in X\}\sub X^n.
$$
Once we know this, we shall find a small (tubular) neighborhood
$U\sub X\hhbr{\Lam}$ of
\begin{equation}\label{dia}
X_\Del\hhbr{n}=X\hhbr{n}\times_{X^n}X^n_\Del
\end{equation}
and a fibration
$$\pi\co U \lra X
$$
whose homotopy type is determined by one of its fiber $\pi\upmo(x)$
and the tangent bundle $TX$. Since the homotopy type of the fiber
$\pi\upmo(x)$ is universal (independent of threefold $X$), the
discrepancy $\delta_{[n]}$ thus only depend on the homotopy type of
$TX$. This will lead to a proof that $\delta_{[n]}$, and thus
$\deg[X\hhbr{n}]\virt$, depends only on a universal expression in
Chern numbers of $X$.

\section{Top diagonal of the Hilbert scheme of $\alpha$--points}

In this section, we shall give a smooth parameterization of the
normal slices to the top diagonal in the Hilbert scheme
$X\hhbr{\alpha}$. Before doing this, we shall comment on the
terminology on stratified spaces and their smooth functions.

\subsection{Stratifications of singular spaces}\label{sec2.3}

Since we primarily are interested in (reduced) quasi-projective
schemes and their open subsets, we shall confine ourselves to their
stratifications and functions.

For quasi-projective $W$, we shall only consider stratifications by
Zariski locally closed smooth subvarieties. It is known that every
quasi-projective scheme admits such stratifications. In case we are
given a finite collection $\cR$ of Zariski closed subsets of $W$, we
can find stratifications $\cS$ of $W$ subordinating to $\cR$ in the
sense that each $R\in\cR$ is a union of strata in $\cS$. To find a
canonical such stratification, we can take the smallest such
stratification\footnote{Two stratifications $\cS_1\leq \cS_2$ if
each strata in $\cS_1$ is the union of strata in $\cS_2$.} among
those subordinating to $\cR$. In case in addition we are given a
morphism of schemes $\pi\mh C\to W$, we can find a stratification
$\cS\pri$ of $C$ and a stratification $\cS$ of $W$ so that $\cS$ is
subordinating to $\cR$ and $\pi\mh C\to W$ is a stratified map.
Again, among all such pairs of stratifications there is one that is
the smallest; we call such pair the standard stratification of $C\to
W$ subordinating to $\cR$.

The collection $\cR$ usually arises from the singular loci of a
sheaf $\cE$ of $\cO_W$--modules. To such sheaf $\cE$ we associate a
collection inductively by letting $R_0=W$ and letting $R_i\sub
R_{i-1}$ be the non-locally free locus of the sheaf of
$\cO_{R_{i-1}}$--modules $\cE\times_{\cO_W}\cO_{R_{i-1}}$. We shall
call $\{R_i\}$ the loci of non-locally freeness of $\cE$.

We next look at the standard stratification of $Y\hhalp$ for a
smooth family over a smooth $T$. Obviously, in case $U\sub T$ is a
Zariski open so that $Y\times_T U=Y_0\times U$, then strata of the
standard stratification of $Y\hhalp\times_TU$ are of the form
$S\times U$ for $S$ strata of $Y_0\hhalp$. In case $Y/T$ is the
total space of a rank three vector bundle, then $Y_0=\At$ and we
know that all strata of $Y_0\hhalp$ are invariant under the symmetry
group of $Y_0$. This proves

\begin{lemm}\label{St1} Let $Y/T$ be a smooth Zariski fiber bundle. Then the
standard stratification of $Y\hhalp$ restricts to the standard
stratification of $(Y_t)\hhalp=(Y\hhalp)_t$.
\end{lemm}

Another property of the standard stratification of $X\hhbr{\alpha}$
is as follows.

\begin{lemm}\label{St2} Let $U_a\sub X$ and $U_b\sub X$ be two analytic open
subset that are biholomorphic to each other, say via $f\co  U_a\lra
U_b$. Then pulling back via $f\sta$ of the universal family $\cI$ of
$X\hhbr{\alpha}$ over $U_b\hhbr{\alpha}\defeq
X\hhbr{\alpha}\times_{X^\Lam}U_a^\Lam$ defines a map
$$\cF\co  U_b\hhbr{\alpha}\lra U_a\hhbr{\alpha}
$$
that is an isomorphism of (standardly) stratified spaces. This
construction also applies to smooth family of quasi-projective
varieties over a smooth base.
\end{lemm}

Here we say a continuous map $f\co  Z_1\to Z_2$ between two stratified
spaces preserves stratifications if every stratum of $Z_1$ is the
preimage of a stratum of $Z_2$. We say $f$ is a smooth isomorphism
of stratified spaces if $f\upmo$ exists and both $f$ and $f\upmo$
are smooth and preserve stratifications of $Z_1$ and
$Z_2$.

\begin{proof} Let $\alpha=(\alpha_1,\cdots,\alpha_k)$. Because
$U_a\hhalp$ is canonically isomorphic to
$$\prod U^{\alpha_i}\times_{U_a^{(\alpha)}} U_a^{[\alpha_i]},
$$
to prove the lemma it suffices to show that for every $n$ the
Hilbert scheme of $n$--points $I_{U_a}(0,n)\cong I_{U_b}(0,n)$
canonically under the map induced by $f$. But this follows from the
universal property of Hilbert schemes.
\end{proof}

\subsection{Smooth functions on stratified spaces}

We now define the notion of smooth functions on a stratified space.

\begin{defi} \label{Def2.1} Let $Z$ be a topological space
equipped with a stratification $\cS$. A smooth function $f\mh
Z\to\CC$ is a continuous function whose restriction to each stratum
is smooth.
\end{defi}

We let $W$ be an analytic space and let $\cE$ be a sheaf of
$\cO_W$--modules. We let $\cS$ be a standard stratification
subordinating to the loci of the non-locally freeness of $\cE$. To
such stratification, we denote by $\cA_\cS$ the sheaf of smooth
functions on $W$, and define the sheaf of smooth sections of $\cE$
be
$$\cA_\cS(\cE)=\cE\otimes_{\cO_W}\cA_\cS.
$$
Because $\cA_\cS$ has partition of unity, whenever
$0\to\cE_2\to\cE\to\cE_1\to 0$ is an exact sequence of sheaves of
$\cO_W$--modules, then both
$$\cA_\cS(\cE_2)\lra\cA_\cS(\cE)\lra\cA_\cS(\cE_1)\lra 0
$$
and
$$\Gamma\bl\cA_\cS(\cE_2)\br\lra\Gamma\bl\cA_\cS(\cE)\br\lra
\Gamma\bl\cA_\cS(\cE_1)\br\lra 0
$$
are exact.

As usual, for any subset $B\sub W$ and sheaf of $\cO_W$--modules
$\cE$, the space of smooth sections of $\cE$ over $B$ is the limit
$$\Gamma(B,\cA_\cS(\cE))=\lim_{B\sub U}\Gamma(U,\cA_\cS(\cE))
$$
taken over all open subsets $U$ containing $B$. Thus each $s$ in
$\Gamma(B,\cA_\cS(\cE)$ is a smooth section defined on some open $U$
containing $B$. Further, in case $\tilde s\in \Gamma(B\pri,
\cA_\cS(\cE))$ is a section over a larger subset $B\pri\supset B$,
we say $s\pri$ extends $s$ if there is an open $U\supset B$ so that
both $s$ and $s\pri$ are defined over $U$ and $s|_{U}\equiv
s\pri|_{U}$. Following this convention, the restriction homomorphism
$$\Gamma(B\pri,\cA_\cS(\cE))\to\Gamma(B,\cA_\cS(\cE))
$$
is surjective for any pair of sets $B\sub B\pri$ with $B$ closed in
$B\pri$,

\subsection{Normal slices to the diagonal in $X\times X$}

We begin with constructing an isomorphism of a tubular neighborhood
$\cU$ of $\Delta(X)\sub X\times X$ with a tubular neighborhood $\cV$
of the zero section $0_X\sub TX$. This map will be a smooth
isomorphism of analytic fiber bundle $\cU/X$ and $\cV/X$ in the
sense that for each $x\in X$ the fiber $\cU_x$ and $\cV_x$ are
isomorphic as analytic spaces. Here $\cU/X$ is via the first
projection $\pr_1\mh \cU\sub X\times X\lra X$.

Since $\cV_x\sub T_xX$, it has a distinguished point $0\in \cV_x$
and $T_0\cV_x\equiv T_xX$ canonically; likewise, $\cU_x$ is a
neighborhood of $x\in X$, which has its distinguished point $x\in
\cU_x$ and isomorphism $T_x\cU_x\equiv T_xX$.

\begin{lemm}\label{lem2.1}
We can find a smooth isomorphism
$$\varphi\co \cU\llra \cV
$$
of a tubular neighborhood $\cU$ of $\Delta(X)\sub X\times X$ and a
tubular neighborhood $\cV$ of $\,0_X\sub TX$, both considered as
vector bundles over $X$, such that
\begin{enumerate}
\item restricting to each fiber $\cU_x$ the map
$\varphi_x=\varphi|_{\cU_x}\mh\cU_x\to \cV_x$ is a biholomorphism,
\item $\varphi_x(x)=0$ and that $d\varphi_x\mh T_x\cU_x\to
T_0\cV_x$ is the identity map.
\end{enumerate}
\end{lemm}

\begin{proof}
First, after identifying $T_xV\cong V$ in the standard way for the
threefold $V=\At$ we can define such $\varphi$ globally:
$$\phi\co  V\times V\lra TV; \quad V\times V\ni (x,v)\longmapsto v-x\in T_xV.
$$
Hence for any open $U_a\sub X$ that admits an open embedding
$U_a\sub V$, the map $\phi$ restricting to
$$\cU\lalp=\phi\upmo(T_{U_a}V)\cap (U_a\times U_a)\sub
X\times X
$$
is a local version of the map required by the lemma; we denote such
$\phi|_{\cU\lalp}$ by
\begin{equation}\label{10.1} \varphi_a\co \cU_a\sub
X\times X\lra TX.
\end{equation}
We next pick an open covering $\{U_a\}$ of $X$ with embedding
$U_a\sub V$ and a partition of unity $\sum\eta_a\equiv 1$
subordinate to the covering $U_a$; we form
$$\cU\defeq\{(x,y)\in X\times X\mid (x,y)\in \cU_a \
\text{whenever}\ x\in \supp(\eta\lalp) \}.
$$
Because $\supp(\eta\lalp)$, which is the closure of $\eta\lalp\ne0$,
is compact and is contained in $U_a$, $\cU$ is open and contains the
diagonal $\Delta(X)\sub X\times X$. Over $\cU$, we define
$$\varphi(x,y)=\sum_{\alpha}\eta_a(x)\varphi_{a}(x,y)\in
TX.
$$
We let $\pi\mh \cU\to X$ be the projection induced by the second
projection $\pr_2\mh X\times X\to X$. Clearly, restricting to each
fiber $\cU_x=\pi\upmo(x)$,
$$\varphi_x(\cdot)=\varphi(x,\cdot)\co  \cU_x\to T_xX
$$
is analytic. Further, because $\varphi_{a}(x,x)=0$ and
$d_y\varphi_{a}(x,y)|_{y=x}=\text{id}$ whenever $\eta_a(x)\ne 0$, we
have $\varphi_x(x)=0$ and $d\varphi_x\mh T_x\cU_x\to T_0\cV_x$ is
the identity map. Therefore, by replacing $\cU$ by a sufficiently
small tubular neighborhood $\cU\pri$ of $\Delta(X)\sub X\times X$
and let $\cV\pri=\varphi(\cU\pri)$, the restriction
$\varphi\pri\defeq\varphi|_{\cU\pri}$ becomes the desired map.
\end{proof}

\subsection{Induced map on Hilbert schemes}

The map $\varphi$ induces a smooth map from the relative Hilbert
scheme $\cU\hhbr{\alpha}$ to $\cV\hhbr{\alpha}$ as fiber bundles
over $X$.

We begin with the definition of $\cU\hhbr{\alpha}$ and
$\cV\hhbr{\alpha}$. First, by letting $Y=TX$ be the vector bundle
over $T=X$, we can apply the previous discussion to form the
relative Hilbert scheme of $\alpha$--points $(TX)\hhbr{\alpha}$. Let
$(TX)\hhbr{\alpha}\to (TX)^\Lam$ be the tautological map; let
$\cV^\Lam$ be the product of $|\Lambda|$--copies of $\cV$ over $X$
indexed by $\Lam$, which is an open subset of $(TX)^\Lam$. We define
$\cV\hhbr{\alpha}$, called the relative Hilbert scheme of
$\alpha$--points in $\cV/X$, be an open subset of $(TX)\hhbr{\alpha}$
defined via the Cartesian product (to the left below)
\begin{equation}
\begin{CD}\label{sq2}
\cV\hhbr{\alpha} @>>> (TX)\hhbr{\alpha}\\
@VVV@VVV\\
\cV^\Lam @>{\jmath_\Lam}>> (TX)^\Lam
\end{CD}\qquad\qquad
\begin{CD}
\cU\hhbr{\alpha} @>>> X\times X\hhbr{\alpha}\\
@VVV@VVV\\
\cU^\Lam @>{\iota_\Lam}>> X\times X^\Lam
\end{CD}
\end{equation}
Similarly we form the relative Hilbert scheme $\cualp$ of
$\alpha$--points of $\cU/X$ as an open subset of $X\times
X\hhbr{\alpha}$ defined by the square (to the right above). Both
$\cualp$ and $\cU\hbr{\alpha}$ are analytic spaces over $X$ via the
first projection of $X\times X\hhalp$ and $X\times X^\Lam$.

Because for each closed $x\in X$, the fiber $\varphi_x\mh
\cU_x\to\cV_x$ is an isomorphism, restricting to fibers over $x$ the
universal families of $\cU\hhbr{\alpha}$ and $\cV\hhbr{\alpha}$
induces a canonical isomorphism $\varphi\hhbr{\alpha}_x\co 
\cU\hhbr{\alpha}_x\lra \cV\hhbr{\alpha}_x$ from
$\cU\hhbr{\alpha}_x=\cU\hhbr{\alpha}\times_X x$ to
$\cV\hhbr{\alpha}_x=\cvalp\times_X x$. We define
\begin{equation}
\varphi\hhbr{\alpha}\co  \cU\hhbr{\alpha}\lra \cV\hhbr{\alpha}
\end{equation}
be $\varphi_x\hhbr{\alpha}$ when restricting to
$\cU\hhbr{\alpha}_x$.

\begin{lemm}
The map $\varphi\hhbr{\alpha}$ is a smooth isomorphism of the
analytic fiber bundles $\cU\hhbr{\alpha}$ and $\cV\hhbr{\alpha}$ as
stratified spaces. Namely, $\varphi\hhbr{\alpha}$ preserves the
standard stratifications of $\cV\hhbr{\alpha}$ and
$\cU\hhbr{\alpha}$ and
$$(\varphi\hhbr{\alpha})\upmo
\bl\cA_{\cV\hhbr{\alpha}}\br=\cA_{\cU\hhbr{\alpha}}.
$$
Further, restricting to fibers over each $x$, $\varphi\hhbr{\alpha}$
induces isomorphism of analytic spaces $\cU_x\hhalp$ and
$\cV_x\hhalp$.
\end{lemm}

\begin{proof}
The proof is a tautology after applying the universal property of
the moduli spaces $\cU\hhbr{\alpha}$ and $\cV\hhbr{\alpha}$. First,
because $\varphi\mh \cU\to\cV$ is a smooth isomorphism that
preserves the analytic structures of the fibers, pulling back the
universal family of $\cV\hhbr{\alpha}$ via $\varphi$ forms a
continuous family of relative $\alpha$--schemes of the fiber bundle
$\cU/X$, thus defines a continuous map
$$\psi\hhbr{\alpha}\co  \cV\hhbr{\alpha}\lra \cU\hhbr{\alpha}.
$$
Now, let $S\sub \cV\hhbr{\alpha}$ be any stratum of the standard
stratification of $\cvalp$. Because its (standard) stratification is
induced from that of $(TX)\hhbr{\alpha}$, the stratum $S$ is the
restriction of a stratum $\tilde S$ of $(TX)\hhbr{\alpha}$. By \fullref{St1}, for an $x\in T$ the intersection $\tilde S\cap
(TX)\hhbr{\alpha}_x$ is a stratum of $(T_xX)\hhbr{\alpha}$. Hence
$S_x=S\cap \cvalp_x$ is a stratum of $\cvalp_x$. Then because
$\psi_x\hhbr{\alpha}\mh \cvalp_x\to\cualp_x$ is an analytic
isomorphism, $S_x\pri=\psi_x\hhbr{\alpha}(S_x)$ must be a stratum of
$\cualp_x$. Similar to the case $TX/X$, $S\pri_x$ is the
intersection with $\cualp_x$ of a single stratum $\tilde S\pri$ of
$X\times X\hhbr{\alpha}$. Let $S\pri=\tilde S\pri\cap\cualp$.
Applying \fullref{St2}, we see immediately that
$$\psi\hhbr{\alpha}(S)=S\pri.
$$
Thus $\psi\hhalp$ preserves the stratifications.

On top of this, because restricting to $S$ the pull back of the
universal family is a smooth family of relative $\alpha$--subschemes,
$\psi\hhbr{\alpha}|_S$ is a smooth map from $S$ to
$\psi\hhbr{\alpha}(S)$. This proves that $\psi\hhbr{\alpha}$ is a
smooth stratified map.

Similarly, applying the same argument to $\varphi\upmo\mh
\cV\to\cU$, we obtain the map $\varphi\hhbr{\alpha}$ that is also a
smooth stratified map. Then the composition
$$\psi\hhbr{\alpha}\circ\varphi\hhbr{\alpha}
\co \cU\hhbr{\alpha}\lra \cU\hhbr{\alpha}
$$
must be the identity map since it is induced by the identity
$\varphi\circ\varphi\upmo=\text{id}$. For the same reason,
$\varphi\hhbr{\alpha}\circ\psi\hhbr{\alpha}=\text{id}$ as well. This
proves that both $\varphi\hhbr{\alpha}$ and $\psi\hhbr{\alpha}$ are
smooth isomorphisms of stratified spaces. Lastly, because
$\psi\hhbr{\alpha}$ preserves fibers over $x$, it defines a smooth
map from $\cV\hhbr{\alpha}$ to $\cU\hhbr{\alpha}$ as smooth fiber
bundles over $X$.

The last statement is true because restricting to fibers over $x\in
X$, both $\varphi\hhbr{\alpha}$ and $\psi\hhbr{\alpha}$ are analytic
isomorphisms. This proves the Lemma.
\end{proof}

\subsection{Top diagonal and its normal slices}

We define the top diagonal $X_\Del\hhbr{\alpha}\sub X\hhbr{\alpha}$
be
$$X_\Del\hhbr{\alpha}=X_\Del^\Lam\times_{X^\Lam}
X\hhbr{\alpha}\qquad \text{where}\quad X_\Del^\Lam=\{(x,\cdots,x)\in
X^\Lam\mid x\in X\}.
$$
Because $X_\Del^\Lam\cong X$, the top diagonal is fibred over $X$:
\begin{equation}\label{20.2}
X_\Del\hhbr{\alpha}\llra X.
\end{equation}
The purpose of this subsection is to find an open neighborhood $U$
of $X_\Del\hhbr{\alpha}\sub X\hhbr{\alpha}$ and a fiber bundle map
$U\to X$ extending the map \eqref{20.2}.

We let $(TX)\hhbr{\alpha}_0$ be the relative Hilbert scheme of
centered $\alpha$--points of the fiber bundle $TX/X$. Because
$\cV\hhbr{\alpha}$ is an open subbundle of $(TX)\hhbr{\alpha}$, we
define the relative Hilbert scheme of centered $\alpha$--points of
$\cV/X$ be
$$\cV\hhbr{\alpha}_0=\cV\hhbr{\alpha}\cap
(TX)\hhbr{\alpha}_0.
$$
There is another way to define this space. Let
$(TX)^\Lam_0=(TX)^\Lam\times_{TX}0_{TX}$ and let $\cV^\Lam_0$ be the
intersection $\cV^\Lam\cap (TX)^\Lam_0$. Then $\cvalp_0$ and
$\cV^\Lam_0$ fit into the following commutative Cartesian squares
\begin{equation}\label{sq3}
\begin{CD} \cvalp_0 @>>>\cV\hhbr{\alpha} @>{\psi\hhbr{\alpha}}>>
\cU\hhbr{\alpha} @>{\jmath\hhbr{\alpha}}>>
X\hhbr{\alpha}\times X@>{\pr_1}>> X\hhbr{\alpha}\\
@VVV@VVV@VVV@VVV@VVV\\
\cV^\Lam_0 @>>>\cV^\Lam @>{\psi_\Lam}>> \cU^\Lam
@>{\jmath_\Lam}>> X^\Lam\times X@>{\pr_1}>> X^\Lam\\
\end{CD}
\end{equation}
in which the $\psi_\Lam$, $\jmath_\Lam$ and $\jmath\hhbr{\alpha}$
are the obvious maps induced by $\psi$ and by $\jmath\mh \cU\to
X\times X$. We let $\Psi\lalp\mh \cvalp_0\to X\hhbr{\alpha}$ be the
composite of the arrows in the top line.

\begin{lemm}\label{le2.6}
After shrinking $\cV$ if necessary, the image set
$$U\hhalp\defeq\Psi\lalp\bl\cV_0\hhbr{\alpha}\br
\sub X\hhbr{\alpha}
$$
is an open neighborhood of $X\hhbr{\alpha}_\Delta\sub
X\hhbr{\alpha}$. Further, after endowing $U$ with the induced
stratification through the inclusion $U\hhalp\sub X\hhbr{\alpha}$,
the induced map
$$\Psi\lalp\co  \cvalp_0\llra U\hhalp
$$
becomes a smooth isomorphism of stratified spaces.
\end{lemm}

\begin{proof}
Because each square in the above diagram is a Cartesian product
square, to prove that $\Psi\lalp\bl\cvalp_0\br$ is open in
$X\hhbr{\alpha}$ it suffices to show that the image of $\cV_0^\Lam$
in $X^\Lam$ under the composite of the bottom line is open in
$X^\Lam$.

We let $h\mh \cV_0^\Lam\to X^\Lam$ be the composite of the bottom
line in the above diagram. We first prove that the differential
$$dh(\xi)\co  T_\xi \cV_0^\Lam\llra T_{h(\xi)} X^\Lam
$$
is an isomorphism at each $\xi$ in the zero section
$0_{(TX)^\Lam}\sub (TX)^\Lam$.

Let $\xi\in 0_{(TX)^\Lam}$ be a point over $x\in X$. Since
$(T_xX)_0^\Lam$ intersects $0_{(TX)^\Lam}$ transversal at $\xi$, the
tangent space $T_{\xi}\cV_0^\Lam$ is a direct sum of
$T_\xi(T_xX)_0^\Lam$ and $T_\xi 0_{(TX)^\Lam}$. Clearly, the images
$d\,h\bl T_\xi(T_xX)_0^\Lam\br$ and $d\,h\bl T_\xi 0_{(TX)^\Lam}\br$
are
$$\{(v_a)_{a\in\cP_\Lam}\in (T_xX)^\Lam\mid \sum v_a=0\} \and
\{(v,\cdots,v)\in (T_xX)^\Lam\mid v\in T_xX\}
$$
respectively. Thus $d\,h$ is an isomorphism at $\xi\in
0_{(TX)^\Lam}\sub \cV_0^\Lam$ and thus $h$ is a diffeomorphism near
$0_{(TX)^\Lam}\sub \cV_0^\Lam$. In particular, if we shrink $\cV$ if
necessary, $h$ becomes a diffeomorphism from $\cV_0^\Lam$ to the
image $h\bl\cV_0^\Lam\br$, an open neighborhood of $X_\Del^\Lam\sub
X^\Lam$.

It follows then that $U\hhalp=\Psi\lalp\bl\cvalp_0\br$ is an open
neighborhood of $X_\Del\hhalp\sub X\hhalp$ and the induced map
$\Psi\lalp$ is continuous and is one--one and onto. We endow
$U\hhalp=\Psi\lalp\bl\cvalp_0\br$ the stratification induced from
that of $X\hhalp$.

We now prove that $\Psi\lalp$ is a smooth isomorphism of stratified
spaces. We first prove that $\Psi\lalp$ preserves the
stratifications. Indeed, because all but the furthest left arrow in
the top line of the above diagram preserves stratifications, to
prove the claim we only need to show that the inclusion
$\cvalp_0\sub\cvalp$ preserves stratifications, which follows from
that the inclusion $(TX)\hhalp_0\to (TX)\hhalp$ preserves
stratifications. For the later, because
$(TX)\hhalp_0=(TX)\hhalp\times_{TX} 0_{(TX)^\Lam}$ with
$\varpi\mh(TX)\hhalp\to TX$ the fiberwise averaging morphism, we
only need to prove that the restriction of $\varpi$ to each stratum
$S\sub (TX)\hhalp$
$$\varpi|_S\co  S\llra TX
$$
is a submersion. But this is clear since the stratum $S$ is induced
by a stratum $S_0$ of $(T_xX)\hhalp$, as shown in \fullref{St1},
and the stratum $S_0$ is invariant under the translation group of
$T_xX$. Therefore $\varpi|_S$ is a submersion.

Once we know that $\Psi\lalp\mh \cvalp_0\to X\hhalp$ preserves
stratifications, the fact that each squares above is a Cartesian
product shows immediately that it is smooth, and its restriction to
each stratum is a diffeomorphism onto its image. Hence, $\Psi\lalp$
is a smooth isomorphism of stratified spaces.
\end{proof}

From the definition, $(TX)\hhalp_0\sub (TX)\hhalp$ is a subscheme,
hence $\cvalp_0\sub\cvalp$ is an analytic subscheme. Therefore,
restricting to fibers $\cualp_{0,x}$ of $\cualp_0$ over $x\in X$ the
map $\Phi$ is analytic.

Before we close this subsection, we shall comment on the partial
equivalence of $\cvalp_0$ and of $U\hhbr{\alpha}$. For this, we
define
$$\cvalp_{0,(\alpha,\beta)}=\cvalp_0\cap
(TX)\hhalp_{0,(\alpha,\beta)} \and
U_{(\alpha,\beta)}\hhalp=U\hhalp\cap X\hhalp\lralpbe.
$$

\begin{lemm}
Let $\Psi\lalp\mh \cvalp_0\to U\hhalp$ be the smooth isomorphism
constructed before. Then
$\Psi\lalp\bl\cvalp_{0,(\alpha,\beta)}\br=U_{(\alpha,\beta)}\hhalp$
and the partial equivalence of $\cvalp_{0}$ and of $U\hhalp$ are
compatible under $\Psi\lalp$.
\end{lemm}

\begin{proof}
This is obvious and will be omitted.
\end{proof}

\subsection{Universal family of Hilbert schemes of centered $\alpha$--points}

We let $Q\to Gr$ be the total space of the universal quotient rank
three vector bundle over the Grassmannian $Gr=Gr(\CC^N,\CC^3)$ and
let $\rho\mh Q\hhalp_0\to Gr$ be the associated relative Hilbert
scheme of centered $\alpha$--points. Let
$$g\co  X\llra Gr
$$
be a smooth map so that $TX\cong g\sta Q$ as smooth vector bundles.
To such $g$, we define the pull back
$$g\sta Q\hhalp_0=Q\hhalp_0\times_{Gr} X=\coprod_{x\in X} \rho\upmo(g(x))
\sub Q\hhalp_0\times X.
$$
Further, the universal family of $Q\hhalp_0$ pulls back under $g$ to
a continuous family of relative centered $\alpha$--points of $TX/X$,
thus defines a continuous map over $X$
\begin{equation}\label{8.2}
g\lalp\co  g\sta Q\hhalp_0\llra (TX)\hhalp_0.
\end{equation}

\begin{lemm}\label{le2.8}
This map is a smooth isomorphism of fiber bundles over $X$. Further,
its restriction to each fiber is an analytic isomorphism and
preserves the partial equivalences of $Q\hhalp_0$ and of
$(TX)_0\hhalp$.
\end{lemm}

\section{Obstruction sheaves}

Our next step is to investigate the obstruction theory of the
relative Hilbert scheme of $\alpha$--points of a smooth family $Y/T$.
Because $Y^\Lam\to Y^{(\alpha)}$ is not \'etale, the defining square
$$\begin{CD}
Y\hhalp @>>> Y^{[\alpha]}\\
@VVV@VVV\\
Y^\Lam @>>> Y^{(\alpha)}
\end{CD}
$$
does not allow us to lift the obstruction theory of $Y^{[\alpha]}$
to that of $Y\hhalp$. Our solution is to take the pull-back of the
obstruction sheaf and the normal cone of $Y^{[\alpha]}$ as the
obstruction sheaf and normal cone of $Y\hhalp$. To force $Y\hhalp\to
Y^{[\alpha]}$ flat, we shall view $Y^\Lam\to Y^{(\alpha)}$, and
hence $Y\hhalp\to Y^{[\alpha]}$, as a morphism between stacks.

We will investigate the obstruction sheaves in this section,
deferring normal cones to the next section. We begin with a brief
account of the relevant extension sheaves.

\subsection{A brief account of extensions sheaves}

Let $Y/T$ be a smooth family of quasi-projective threefolds over a
smooth variety $T$; let $S$ be any scheme over $T$ and $Z$ be a flat
$S$--family of 0--subschemes in $Y_S=Y\times_TS$; and let $p_S\mh
Y_S\to S$ be the projection. Our first goal is to prove a canonical
isomorphism relating the traceless part of the relative extension
sheaf of the ideal sheaf $\cI_Z$ of $Z\sub Y_S$ with the high direct
image sheaf of an extension sheaf.

\begin{lemm}\label{lem3.1}
Let the notation be as before, then we have canonical isomorphism
$$\ext^2_{p_S}\bl \cI_Z,\cI_Z\br_0\cong p_{S\ast}
\ext^3_{\cO_{Y_S}}(\cO_Z,\cI_Z).
$$
\end{lemm}

\begin{proof}
We begin with the defining exact sequence
$$0\lra\cI_Z\lra\cO_{Y_S}\lra\cO_Z\lra 0
$$
and its associated long exact sequence of relative extension sheaves
$$\lra \ext^2_{p_S}(\cO_{Y_S},\cI_Z)\mapright{\phi_2}
\ext^2_{p_S}(\cI_Z,\cI_Z)\lra \ext^3_{p_S}(\cO_Z,\cI_Z)\lra
\qquad\qquad
$$
\begin{equation}
\label{3.1} \qquad\qquad\qquad\lra\ext^3_{p_S}(\cO_{Y_S},\cI_Z)
\mapright{\phi_3}\ext^3_{p_S}(\cI_Z,\cI_Z)\lra0.
\end{equation}
Because $\cO_{Y_S}$ is locally free and $Z$ is a flat family of
zero-subschemes $Y_S$,
$$\ext^k_{p_S}(\cO_{Y_S},\cI_Z)=
R^kp_{S\ast}\cI_Z=R^kp_{S\ast}\cO_{Y_S},\quad k\geq 2.
$$
On the other hand, it follows from the definition that
$R^kp_{S\ast}\cO_{Y_S}$ is a subsheaf of $\ext^k_{p_S}(\cI_Z,\cI_Z)$
and the composite
\begin{equation}
\label{3.2} R^kp_{S\ast}\cO_{Y_S}
\mapright{\sub}\ext^k_{p_S}(\cI_Z,\cI_Z)
\mapright{\text{tr}}R^kp_{S\ast}\cO_{Y_S}
\end{equation}
is multiplying by $1$, which is the rank of $\cI_Z$. Thus
$R^kp_{S\ast}\cO_{Y_S}$ is the trace part of
$\ext^k_{p_S}(\cI_Z,\cI_Z)$, and the traceless part
\begin{equation}\label{3.3}
\ext^2_{p_S}(\cI_Z,\cI_Z)_0\cong\ext^{3}_{p_{S}}(\cO_Z,\cI_Z)
\end{equation}
canonically.

To complete the proof of the lemma, we shall apply the local to
global spectral sequence. First because the extension sheaf
$\ext^i(\cO_Z,\cI_Z)$ is zero away from $Z$ and $Z\to S$ has
relative dimension $0$,
$$R^kp_{S\ast}\ext^{3-k}(\cO_Z,\cI_Z)=0
$$
except when $k=0$. Hence the spectral sequence
$$R^{\bullet}p_{S\ast}\ext^{\bullet}(\cO_Z,\cI_Z)
\Rightarrow \ext^{\bullet}_{p_S}(\cO_Z,\cI_Z)
$$
degenerates to
$$\ext^3_{p_S}\bl \cO_Z,\cI_Z\br_0\cong p_{S\ast}
\ext^3(\cO_Z,\cI_Z).
$$
This proves the Lemma.
\end{proof}

The relative extension sheaves satisfy the base change property. Let
$\rho\mh S\pri\to S$ be a morphism and $Z\pri=Z\times_S S\pri$ be
the subscheme of $Y_{S\pri}=Y_S\times_SS\pri$, then
\begin{equation}\label{base}
\ext^2_{p_{S\pri}}\bl \cI_{Z\pri},\cI_{Z\pri}\br_0\cong \rho\sta
\ext^2_{p_S}\bl \cI_Z,\cI_Z\br_0
\end{equation}
and such isomorphisms are compatible with $S^{\prime\prime}\to
S\pri\to S$.

To prove the base change property, we will apply the Lemma and prove
the base change property of
$$p_{S\ast}\ext^3(\cO_Z,\cI_Z).
$$
Because $\dim Y/T=3$ and $\cO_Z$ is flat over $S$, the later admits
a length three locally free resolution as a sheaf of
$\cO_{Y_S}$--modules. Then because tensoring this resolution by
$\cO_{S\pri}$ provides a locally free resolution of $\cO_{Z\pri}$,
applying the definition of extension sheaves we readily see that
\begin{equation}\label{base2}
\ext^3_{\cO_{Y_{S\pri}}}(\cO_{Z\pri},\cI_{Z\pri})\cong
\ext^3_{\cO_{Y_{S}}}(\cO_{Z},\cI_{Z})\otimes_{\cO_{S}}\cO_{{S\pri}}.
\end{equation}
Because the first sheaf if a finite sheaf of $\cO_{S\pri}$--modules,
$$p_{S\pri\ast}\ext^3_{\cO_{Y_{S\pri}}}(\cO_{Z\pri},\cI_{Z\pri})
=\ext^3_{\cO_{Y_{S\pri}}}(\cO_{Z\pri},\cI_{Z\pri})
$$
when the later is viewed as sheaf of $\cO_{S\pri}$--modules. Thus by
viewing the other $\ext^3$ sheaf as a sheaf of $\cO_S$--modules, the
identity is exactly the base change property \eqref{base}.

\begin{coro}\label{cor3.2}
Let $Z\sub Y_S$ be as in the previous lemma. Suppose further that
$Z=Z_1\cup Z_2$ is a union of two disjoint flat $S$--families of
0--subschemes. Then
$$\ext^2_{p_S}\bl \cI_Z,\cI_Z\br_0\cong
\ext^2_{p_S}\bl \cI_{Z_1},\cI_{Z_1}\br_0\oplus \ext^2_{p_S}\bl
\cI_{Z_2},\cI_{Z_2}\br_0
$$
canonically.
\end{coro}

\begin{proof}
This is true because
$$\ext^2_{p_S}\bl \cI_Z,\cI_Z\br_0\cong
p_{S\ast} \ext^3(\cO_Z,\cI_Z)\cong \bigoplus_{i=1}^2 p_{S\ast}
\ext^3(\cO_{Z_i},\cI_{Z_i}).\proved
$$
\end{proof}

\subsection{Obstruction sheaves of Hilbert scheme of points}

For a smooth family of quasi-projective threefolds over a smooth
base $T$, the relative Hilbert scheme $I_{Y/T}(0,n)$ is a
quasi-project fine moduli scheme with universal family
$$\cZ\sub Y\times_T I_{Y/T}(0,n);
$$
its obstruction theory as moduli of stable sheaves with fixed
determinants is perfect with obstruction sheaf the traceless part of
the relative extension sheaf under the second projection of
$Y\times_T I_{Y/T}(0,n)$:
$$\Ob_{I_{Y/T}(0,n)}=\ext^2_{\pr_2}(\cI_{\cZ},\cI_{\cZ})_0.
$$
For $Y\hbr{n}$, which is $I_{Y/T}(0,n)$ endowed with the reduced
scheme structure, we shall take the pull back of the obstruction
sheaf of $I_{Y/T}(0,n)$ under the inclusion $Y\hbr{n}\to
I_{Y/T}(0,n)$ as its obstruction sheaf. By the base change property
just proved, it can also be expressed as the traceless part of the
relative extension sheaf of the universal ideal sheaf of $Y\hbr{n}$.

Now let $\alpha=(\alpha_1,\cdots,\alpha_k)$ be any element in
$\cP\llam$ as before; let $(\cZ_1,\cdots,\cZ_k)$ be the universal
family of $Y\hbr{\alpha}$, each is a flat family of length
$|\alpha_i|$ zero-subschemes in $Y/T$ over $Y\hbr{\alpha}$; we let
$\cI_{\cZ_i}$ be the ideal sheaf of $\cZ_i\sub Y\times_T
Y\hbr{\alpha}$. Because $Y\hbr{\alpha}=\prod_T Y\hbr{\alpha_i}$, we
define the obstruction sheaf of $Y\hbr{\alpha}$ be the direct sum of
the pull back of that of $Y\hbr{\alpha_i}$, which by \fullref{cor3.2} is of the form
$$\Ob_Y^{[\alpha]}=\bigoplus_{i=1}^k \ext^2_{\pr_2}\bl
\cI_{\cZ_i},\cI_{\cZ_i}\br_0.
$$
As to $Y\hhalp$, we shall take the pull-back $\phi\lalp\sta
\Ob\hbr{\alpha}$ (of the tautological $\phi\lalp\mh
Y\hhbr{\alpha}\to Y\hbr{\alpha}$) as its obstruction sheaf. For the
same reason, they can be defined using relative extension sheaves.
Let $(\cW_1,\cdots, \cW_k)$ with $\cW_i\sub Y\times_T Y\hhalp$ be
part of the universal family of $Y\hhbr{\alpha}$; each $\cW_i$ is
the pull back of $\cZ_i$ under $Y\hhalp\to Y^{[\alpha]}$; let
$\cI_{\cW_i}$ be the ideal sheaf of $\cW_i\sub Y\times_T
Y\hhbr{\alpha}$. The obstruction sheaf of $Y\hhalp$ is
$$\Ob_Y\hhbr{\alpha}\equiv \bigoplus_{i=1}^k \ext^2_{\pr_2}(\cI_{\cW_i},
\cI_{\cW_i})_0.
$$

\subsection{Comparing obstruction sheaves under
equivalences}\label{sec3.3}

Our next task is to compare the sheaves $\Ob_Y\hhbr{\alpha}$ with
$\Ob_Y\hhbr{\beta}$ over the partial equivalence
$Y\hhalp\lralpbe\cong Y\hhbe\lrbealp$.

\begin{lemm}\label{she-iso}
Under the partial equivalence
$\hhy{\alpha}\lralpbe\cong\hhy{\beta}_{(\beta,\alpha)}$, the
restriction to $\hhy{\alpha}\lralpbe$ of the obstruction sheaf
$\Ob_Y\hhbr{\alpha}$ is canonically isomorphic to the restriction to
$\hhy{\beta}_{(\beta,\alpha)}$ of $\Ob_Y\hhbr{\beta}$.
\end{lemm}

\begin{proof}
We first show that the Lemma can be reduced to the case where
$\alpha\geq\beta$. Indeed, because $\hhy{\alpha}\lralpbe\cong
\hhy{\beta}_{(\beta,\alpha)}$ is induced by
$\hhy{\alpha}\lralpbe\sub \hhy{\alpha}_{(\alpha,\awb)}\cong
\hhy{\awb}_{(\awb,\alpha)}$, should the lemma hold for
$\alpha\geq\beta$, we would have
$$\Ob_Y\hhbr{\alpha}|_{\hhy{\alpha}_{(\alpha,\awb)}}\cong
\Ob_Y\hhbr{\awb}|_{\hhy{\awb}_{(\awb,\alpha)}},
$$
which would imply
$$\Ob_Y\hhbr{\alpha}|_{\hhy{\alpha}_{(\alpha,\beta)}}\cong
\Ob_Y\hhbr{\awb}|_{\hhy{\alpha}_{(\alpha,\beta)}=\hhy{\beta}_{(\beta,\alpha)}}
\cong\Ob_Y\hhbr{\beta}|_{\hhy{\beta}_{(\beta,\alpha)}}.
$$
For the case $\alpha\geq\beta$, by induction we only need to
consider the case where
$$\alpha=(\alpha_1,\cdots,\alpha_k)\geq
\beta=(\beta_1,\cdots,\beta_{k+1}).
$$
For simplicity we assume $\alpha_i=\beta_i$ for $i<k$ and
$\alpha_k=\beta_k\cup\beta_{k+1}$.

Now let $(\cW_i,\varphi_i)_{1\leq i\leq k}$, with $\cW_i\sub
Y\times_T Y\hhalp\lralpbe$ be the universal family of $\hhy{\alpha}$
over $Y\hhalp\lralpbe$; let
$(\tilde\cW_i,\tilde\varphi_i)_{i=1}^{k+1}$ be the universal family
of $\hhy{\beta}$ over $Y\hhbe\lrbealp$. Because $\alpha_i=\beta_i$
for $i<k$, $\cW_i\cong\tilde\cW_i$ after identifying
$Y\hhalp\lralpbe\cong Y\hhbe\lrbealp$. On the other hand, by
\fullref{cor1.2} the supports of $\tilde\cW_k$ and
$\tilde\cW_{k+1}$ are disjoint closed subsets of $Y\times_T
\hhy{\beta}_{(\beta,\alpha)}$; and their union form a flat family of
zero-subschemes satisfying $\tilde\cW_k\cup\tilde\cW_{k+1}=\cW_k$.

As for the obstruction sheaves $\Ob_Y\hhbr{\alpha}$ and
$\Ob_Y\hhbr{\beta}$, by definition,
\begin{gather*}
\Ob_Y\hhbr{\alpha}|_{Y\hhalp\lralpbe}
=\bigoplus_{i=1}^k\ext^2_{\pr_2}\bl \cI_{\cW_i},\cI_{\cW_i}\br_0
\\
\tag*{\hbox{and}}
\Ob_Y\hhbr{\beta}|_{Y\hhbe\lrbealp}=\bigoplus_{i=1}^{k+1}\ext^2_{\pr_2}\bl
\cI_{\tilde\cW_i},\cI_{\tilde\cW_i}\br_0.
\end{gather*}
Hence to prove the lemma we only need to check that
$$\ext^2_{\pr_2}\bl\cI_{\cW_k},\cI_{\cW_k}\br_0
\cong \ext^2_{\pr_2}\bl \cI_{\tilde\cW_k},\cI_{\tilde\cW_k}\br_0
\oplus \ext^2_{\pr_2}\bl
\cI_{\tilde\cW_{k+1}},\cI_{\tilde\cW_{k+1}}\br_0.
$$
But this is exactly what was proved in the last subsection.
\end{proof}

\subsection{Obstruction sheaves under smooth isomorphism}

We next move to the obstruction sheaves of various moduli spaces of
interests. We continue to denote by $X$ a smooth complex threefold.

Let $\cI_{\cZ_i}$ be the the ideal sheaves of
$\cZ_1,\cdots,\cZ_k\sub \cU_0\hhbr{\alpha}\times_X\cU$, which are
part of the universal family of $X\hhbr{\alpha}$; let $\pi_1$ be the
first projection of $X\hhbr{\alpha}\times X$. Then the obstruction
sheaf of $X\hhalp$ is
\begin{equation}\label{e}
\Ob_X\hhbr{\alpha}=\bigoplus_{i=1}^k
\ext^2_{\pi_1}\bl\cI_{\cZ_i},\cI_{\cZ_i}\br_0 \cong\bigoplus_{i=1}^k
\pi_{1\ast}\ext^2\bl \cO_{\cZ_i},\cI_{\cZ_i}\br
\end{equation}
Similarly, the obstruction sheaf of $\cvalp_0$ is defined to be the
relative extension sheaf
\begin{equation}\label{ex}
\Ob_{\cV_0}\hhbr{\alpha}=\bigoplus_{i=1}^k
\ext^2_{\pi_1}\bl\cI_{\cW_i},\cI_{\cW_i}\br_0 \cong\bigoplus_{i=1}^k
\pi_{1\ast}\ext^2\bl \cO_{\cW_i},\cI_{\cW_i}\br
\end{equation}
in which $\cI_{\cW_i}$ are ideal sheaves of the subschemes
$\cW_1\cdots,\cW_k\sub \cvalp_0\times_X \cV$ that are part of the
universal family of $\cvalp_0$; that $\pi_1$ is the first projection
of $\cvalp_0\times_X \cV$.

These obstruction sheaves are related in the obvious way. First, we
let
$$\cA_1=\cA_{X\hhalp\times X/X\hhalp}
\and \cA_2=\cA_{\cvalp_0\times_X\cV/\cvalp_0}
$$
respectively be the sheaf of smooth functions on $X\hhalp\times X$
and ${\cvalp_0\times_X\cV}$ that are analytic along fibers of
$X\hhalp\times X/X\hhalp$ and ${\cvalp_0\times_X\cV/\cvalp_0}$.
Because $\Psi\lalp$ of \fullref{le2.6} is induced by the universal
family and because the tautological map
$$h=(\Psi\lalp,p_2\circ\psi)\co \cvalp_0\times_X \cV \lra
X\hhalp \times X,
$$
where $p_2\circ \psi\mh \cV\to X$ is the composite of $\psi\mh
\cV\to X\times X$ with the second projection of $X\times X$, maps
fibers to fibers and is analytic along fibers,
\begin{equation}\label{ext}
h\sta\bl \cI_{\cZ_{i}}\otimes \cA_1\br\cong \cI_{\cW_{i}}\otimes
\cA_2.
\end{equation}
Here the tensor products are over the sheaves $\cO_{X\hhalp\times
X}$ and $\cO_{\cvalp_0\times_X\cV}$. On the other hand, by the
property of relative extension sheaf, the tensor product
\begin{eqnarray*}
\Ob_{\cV_0}\hhbr{\alpha}\otimes\cA_{\cvalp_0}&\cong&
\bigoplus_{i=1}^k
\pi_{1\ast}\Bl\ext^3_{\cA_2}(\cI_{\cW_i},\cI_{\cW_i})\otimes
\cA_2\Br\\
&\cong& \bigoplus_{i=1}^k \pi_{1\ast}\ext^3_{\cA_2}\bl
\cI_{\cW_i}\otimes \cA_2,\cI_{\cW_i}\otimes \cA_2\br.
\end{eqnarray*}
Here the extension sheaf $\ext^3_{\cA_2}(\cdot,\cdot)$ is defined
using locally free resolution of locally free sheaves of
$\cA_2$--modules.

Because of the isomorphism \eqref{ext}, the last term in the above
isomorphisms is canonically isomorphic to
\begin{eqnarray*}
&\cong& \bigoplus_{i=1}^k \pi_{1\ast}\Bl\ext^3_{\cA_2}\bl
h\sta(\cI_{\cZ_i}\otimes\cA_1),h\sta(\cI_{\cZ_i}\otimes
\cA_1)\br\Br\\
&\cong& \bigoplus_{i=1}^k \pi_{1\ast}\Bl\ext^3_{\cA_1}\bl
\cO_{\cZ_i}\otimes \cA_1,\cI_{\cZ_i}\otimes \cA_1\Br\ \cong\
\Psi\lalp\sta \Bl\Ob_{X}\hhbr{\alpha}\otimes\cA_{X\hhalp}\Br.
\end{eqnarray*}
This proves that
\begin{equation}\label{9.3}
\Psi\lalp\sta \Bl\Ob_{X}\hhbr{\alpha}\otimes\cA_{X\hhalp}\Br \cong
\Ob_{\cV_0}\hhbr{\alpha}\otimes\cA_{\cvalp_0}.
\end{equation}

\section{Representing Virtual Cycles}

We will begin this section with a quick review of the construction
of the virtual cycle of a scheme with perfect obstruction theory.
For more details of this construction, please consult
Behrend and Fantechi \cite{Behr-Fant} or Li and Tian \cite{Li-Tian}.

\subsection{Virtual cycles via Gysin map}

For the moment, we assume $W$ is a quasi-project scheme with a
perfect obstruction theory and obstruction sheaf $\Ob$. Following
\cite{Behr-Fant,Li-Tian}, the virtual cycle of $W$ is constructed
via
\begin{enumerate}
\item finding a vector bundle $V$ on $W$ and a surjective sheaf
homomorphism\footnote{In this paper, whenever we use a Roman
alphabet, say $V$, to denote a vector bundle, we will use its
counter part $\cV$ to denote the sheaf of regular sections
$\cO_W(V)$.} $\cV\to \Ob$; then the obstruction theory of $W$
provides us a unique cone cycle\footnote{A cycle is a finite union
$\sum m_iD_i$ of subvarieties $D_i$ with integer coefficients $m_i$;
it is a cone cycle if all $D_i$ are cones in $V$.} $C \sub V$ of
codimension $\rk V$;
\item defining the virtual cycle $[W]\virt=0_{V}\sta[C]$ via the Gysin
homomorphism
$$\smash{0_V\sta\mh H\lsta^{BM}(V,\ZZ)\to H\lsta(W,\ZZ)}.$$
Since every subvariety of $V$ defines a class in the Borel--Moore homology
group, $0_V\sta[C]$ is well-defined.
\end{enumerate}
Before we move on, a few comments are in order.

Usually, the Gysin map is defined as a homomorphism between Chow
groups. For us, we shall use the Borel--Moore homology group and use
intersecting with smooth sections to define this homomorphism.

The cone cycle $C\sub V$ is unique in the following sense. To each
closed $w\in W$, we let $\hat w$ be the formal completion of $W$
along $w$ and fix an embedding
$$\hat w\sub T\hat\defeq\spec \kk[\![
T_w W\dual]\!]
$$
that is consistent with the tangent space at their only closed
points. We let $O=\Ob\otimes_{\cO_W}\kk(w)$ be the obstruction space
to deforming $w$ in $W$. Then the Kuranishi map of the obstruction
theory provides a canonical embedding of the normal cone $C_{\hat
w}\hat T$ to $\hat w$ in $\hat T$:
$$C_{\hat w}\hat T\sub \hat w\times O,
$$
where the later is viewed as a vector bundle over $\hat w$. The
uniqueness of $C$ asserts that there is a vector bundle homomorphism
$$\eta_w\co  V\times_W \hat w\llra O\times\hat w
$$
extending the homomorphism $V|_w\to O$ induced by $\cV\to\Ob$ such
that
\begin{equation}\label{cone1}
\eta\sta C_{\hat w}\hat T=C\cap (V\times_W\hat w),
\end{equation}
Lastly, the resulting cycle $0_V\sta[C]\in A\lsta W$ is independent
of the choice of $\cV\to\Ob$.

As we will see in the later part of this paper, it will be useful to
eliminate the dependence of constructing $[W]\virt$ on the choice of
$\cV\to\Ob$. To achieve this, we will use smooth sections of $\Ob$
to define the cycle $[W]\virt$.

Before we do that, we shall first recall the notion of
pseudo-cycles.

\subsection{Pseudo-Cycles}

In this work, we shall use pseudo-cycles to represent homology
classes in a stratifiable space.

Let $\Theta\sub W$ be a triangulable closed subset of a stratified
space $W$. We shall fix a Riemannian metric\footnote{We can embed
$W$ in a smooth space and use the induced Riemannian metric on the
ambient space.} on $W$ and denote by $\Theta_\eps$ the
$\eps$--tubular neighborhood $\Theta$ in $W$.

\begin{defi} A $(\Theta,\eps)$--relative $d$--dimensional pseudo-chain is a
pair $(f,\Sigma)$ of a smooth, oriented $d$--dimensional manifold
with smooth boundary $\partial\Sigma$ and a continuous map $f\mh
\Sigma\to W$ such that $f$ is smooth over
$\Sigma-f\upmo(\Theta_\eps)$.
\end{defi}

We denote the $\ZZ$--linear span of all such pseudo-chains by
$PCh_d(W)_\Theta$ with $\eps$--implicitly understood.

We define an equivalence relation $\cong$ on $PCh_d(W)_\Theta$ as
follows. For notation brevity, we shall use $[f]\pc$ to denote the
pseudo-cycle $(f,\Sigma)$ with $\Sigma$ implicitly understood. We
call two such chains $[f_1]\pc$ and $[f_1]\pc$ equivalent if there
are open subsets $A_i\sub \Sigma_i$ and orientation preserving
diffeomorphism $\phi\mh A_1\to A_2$ so that $\Sigma_i-A_i\sub
f_i\upmo(\Theta_\eps)$ for $i=1$ and $2$, and $f_2\circ
\phi|_{A_1}=f_1|_{A_1}$; we define $[f_1\cup f_2]\pc\cong
[f_1]\pc+[f_2]\pc$ with the union $[f_1\cup f_2]\pc$ be the map
$\Sigma_1\cup\Sigma_2\to W$ induced by $f_1$ and $f_2$; we define
$[-f]\pc\cong-[f]\pc$ with $[-f]\pc$ be $f\mh \Sigma^-\to W$ and
$\Sigma^-$ the space $\Sigma$ with the opposite orientation. The
equivalence relation $\cong$ generates an ideal in the Abelian group
$PCh_d(W)_\Theta$.

There is an obvious boundary homomorphism
$$\partial\co  PCh_d(W)_\Theta\llra PCh_{d-1}(W)_\Theta
$$
that sends any $[f]\pc$ to $[\partial f]\pc$ with $\partial f$ is
the restriction of $f$ to $\partial \Sigma$. The kernel of its
induced homomorphism is defined to be the space of pseudo-cycles
relative to $\Theta$:
$$PCy_d(E)_\Theta\defeq \ker\{\partial\co PCh_d(W)_\Theta\lra
PCh_{d-1}(W)_\Theta/\cong\}.
$$
The proof of the following lemma is standard.

\begin{lemm}\label{Kr} Suppose $\dim_{\RR}\Theta\leq d-2$ and
suppose the $\eps$--tubular neighborhood $\Theta_\eps$ is a
deformation retract to $\Theta$. Then every $(\Theta,\eps)$--relative
pseudo-cycle $[f]\pc\in Pcy_d(W)_\Theta$ defines canonically a
homology class $[f]\in H_d(W,\ZZ)$.
\end{lemm}

One version of pseudo-cycle we shall repeatedly use is the
following:

\begin{rema}\label{pc}
Any pair $(B,\Theta)$ of a closed subset $B\sub W$ and a
stratifiable closed $\Theta\sub W$ such that $B-\Theta$ is a smooth,
oriented $d$--dimensional manifold and $\dim_\RR \Theta\leq d-2$ is a
$(\Theta,\eps)$--relative pseudo-cycle for all $\eps>0$.
\end{rema}

This can be seen as follows. For any $\eps>0$, we pick an open
$O\sub B\cap \Theta_{\eps/2}$ so that $\Sigma=B-O$ is a smooth
manifold with smooth boundary. Then the identity map $f\mh \Sigma\to
B\sub W$ is a $(\Theta,\eps)$--relative pseudo-cycle.

In case $W$ is stratifiable, which is the case when $W$ is an open
subset of a quasi-projective scheme, $\Theta_\eps$ deformation
retract to $\Theta$ for all sufficiently small $\eps$. Because
$\dim_\RR \Theta\leq d-2$, for all sufficiently small $\eps$ these
$(\Theta,\eps)$--relative pseudo cycles all represent the same
homology class in $H_d(W,\ZZ)$. Because of this, in the future we
will call $B$ a $\Theta$--relative pseudo-cycle and will denote by
$[B]\in H_d(W,\ZZ)$ the resulting homology class. In case $\Theta$
is the singular locus of $B$, we shall call $B$ a pseudo-cycle
directly with $\Theta=B_{sing}$ implicitly understood.

\subsection{Cycle representatives via smooth sections of
sheaves}\label{rep4.4}

We now investigate how to intersect with smooth sections to define
the Gysin homomorphism of a cone cycle in a vector bundle over $W$.

Let $C$ be a pure $\CC$--dimension $r$ algebraic cone cycle in a
vector bundle $\pi\mh V\to W$;
$$C=\sum m_i C_i
$$
its irreducible components decomposition. We pick a stratification
$\cS$ of $W$ and a stratification $\cS\pri$ of $\cup_i C_i$ so that
the induced map $\cup_i C_i\to W$ is a stratified map. Also, to each
stratum $S\pri\in\cS\pri$ with $S=\pi(S\pri)$ its image stratum, we
shall denote by $V_{S}$ the restriction to $S$ of $V$; for a smooth
section $s\in \cA_\cS(V)$, we denote its graph by $\Gamma_s\sub V$.

\begin{defi}\label{good}
A smooth section $s\in \cA_\cS(V)$ is said to intersect
transversally with $C$ if every stratum $S\pri\sub \cup_i C_i$,
which is a smooth subset in $V_{S}$, intersects transversally with
$\Gamma_s\cap V_{S}$ inside $V_{S}$.
\end{defi}

By embedding $W$ in a projective space and extending $V$, we can
apply the standard Sard's transversality theorem to conclude that
the set of smooth sections that intersect transversally with $C$ is
dense in the space of all smooth sections.

Let $s$ be a section that intersects transversally with $C$. We
claim that the intersection $\Gamma_s\cap C$ is a pseudo-cycle.
Indeed, let $S_i\sub C_i$ be the open stratum of $C_i$; then
$S_i\sub C_i$ is smooth, open and Zariski dense; therefore,
$\dim_\RR C_i-S_i\leq \dim_\RR C_i-2$. Suppose $\dim C_i=d$ (for all
$i$) and $\rk V=r$. Then $\Gamma_s\cap S_i$ is a smooth, oriented
manifold of real dimension $2d-2r$; its complement
$\Gamma_s\cap(C_i-S_i)$ has real dimension at most $2d-2r-2$.
According to \fullref{pc}, by taking
$$\Theta=\cup_i\Gamma_s\cap(C_i-S_i),
$$
the set $\Gamma_s\cap C_i$ becomes a $2d-2r$--dimensional
pseudo-cycle (relative to $\Theta$). We denote this pseudo-cycle by
$(\Gamma_s\cap C_i)\pc$ and denote its image pseudo-cycle under the
projection $\pi\mh V\to W$ by $\pi\lsta(\Gamma_s\cap C_i)\pc$.
Finally, by linearity,
$$(\Gamma_s\cap C)\pc=\sum_i(\Gamma_s\cap C_i)\pc\and
\pi\lsta(\Gamma_s\cap C)\pc=\sum_i \pi\lsta(\Gamma_s\cap C_i)\pc.
$$
It is immediate to check that its associated homology class given by
\fullref{Kr} is the Gysin homomorphism image of $[C]$
$$0_V\sta[C]=[\pi\lsta(\Gamma_s\cap C)\pc]
\in H_{2d-2r}(W,\ZZ).
$$
We now apply this technique to construct the virtual cycle
$[W]\virt$, assuming that $W$ is quasi-projective with a perfect
obstruction theory and the obstruction sheaf $\Ob$. As we mentioned,
we first pick a locally free sheaf $\cV$ and a quotient sheaf
homomorphism $\cV\to\Ob$, thus obtaining a cone $C\sub V$ given by
the obstruction theory of $W$. We then pick the standard pair of
stratification $\cS\pri$ of $C$ and $\cS$ of $W$ that respects the
morphism $C\to W$ and the loci of non-locally freeness of $\Ob$.
Then as was argued, we can find a section $s\in \cA_\cS(V)$ that
intersects transversally with $C$. The virtual fundamental cycle
$[W]\virt$, which is $0_V\sta[C]$, is the homology class given by
the pseudo-cycle
$$[W]\virt=0_V\sta[C]=[\pi\lsta(C\cap\Gamma_s)\pc]\in H\lsta(W,\ZZ).
$$
We next show that such pseudo-cycle can be constructed using the
image section $\xi\in \cA_\cS(\Ob)$ of $s$ under the homomorphism
$\cA_\cS(V)\to \cA_\cS(\Ob)$. To this end, we need first to recover
the pseudo-cycle $\pi\lsta(\Gamma_s\cap C)\pc$ using the smooth
section $\xi\in \cA_\cS(\Ob)$. Suppose $\cV\pri\to\Ob$ is a
surjective homomorphism and $t\in \cA_\cS(V\pri)$ is a lift of $\xi$
under $\cA_\cS(\cV\pri)\to\cA_\cS(\Ob)$, which exists by the exact
sequence at the end of  subsection 2.2. We claim that $t$ intersects
transversally with the virtual normal cone $C\pri\sub V\pri$ and
\begin{equation}\label{22}\pi\lsta(\Gamma_s\cap C)\pc
=\pi\pri\lsta(\Gamma_t\cap C\pri)\pc
\end{equation}
as pseudo-cycles in $W$.

We first prove the case where $\cV\pri=\cV$ and $t\in \cA_\cS(V)$ is
another lifting of $\xi$. We first show that as sets
$\pi(\Gamma_s\cap C)=\pi\pri(\Gamma_t\cap C\pri)$. For this, we
consider the section $s-t\in \cA_\cS(V)$ and its induced fiberwise
translation
$$\ell_{s-t}\co  V\to V;\qquad x\in V_w\longmapsto x+s(w)-t(w)\in V_w.
$$
Because of \cite{Li-Tian}, the fiber $C\cap V_w$ of the cone over
$w$ is translation invariant under vectors in
$$K_w=\ker\{V_w\to
\Ob|_w\}.
$$
Therefore $\ell_{s-t}\mh V\to V$ maps $C$ to $C$. But on the other
hand, $\ell_{s-t}$ maps $\Gamma_t$ to $\Gamma_s$, hence
$\ell_{s-t}(\Gamma_t\cap C)=\Gamma_s\cap C$, which proves
\eqref{22}.

It remains to show that $t$ intersects $C$ transversally and the
induced orientation on $\pi(\Gamma_s\cap C)$ coincides with that of
$\pi(\Gamma_t\cap C)$. First, because the stratification $\cS$
respects the non-locally freeness of $\Ob$, the restriction
$\Ob\otimes_{\cO_W}\cO_{S}$ is locally free. Hence the collection
$\{K_w\mid w\in S\}$ forms a subbundle of $V_{S}$. By the
translation invariance of $C\cap V_w$ under $K_w$ and by the
minimality of the stratification $\cS\pri$, $S\pri\cap V_w$ is
invariant under the translations by vectors in $K_w$; hence
$S\pri=\ell_{s-t}(S\pri)$.

On the other hand, since the map $\ell_{s-t}\mh V_S\to V_S$ is a
smooth diffeomorphism, $\Gamma_t$ intersects transversally with
$S\pri$ if and only if its image $\ell_{s-t}(\Gamma_t)=\Gamma_s$
intersects transversally with $\ell_{s-t}(S\pri)=S\pri$, which is
true by our choice of $s$. Hence $\Gamma_t$ intersects transversally
with $S\pri$. In particular, $\ell_{s-t}(\Gamma_t\cap
C)=\Gamma_s\cap C$ and hence $\pi\lsta(\Gamma_t\cap
C)\pc=\pi\lsta(\Gamma_s\cap C)\pc$.

Next we consider the general situation where $\cV\pri\to\Ob$ is an
arbitrary quotient sheaf homomorphism by a locally free sheaf. We
claim that by picking a lifting $t\in \cA_\cS(\cV\pri)$ of $\xi\in
\cA_\cS(\Ob)$, we obtain the identical cycle representatives of
$[W]\virt$. Indeed, by our previous discussion, we only need to
consider the case that $\cV$ is a quotient sheaf of $\cV\pri$ and
$\cV\pri\to\Ob$ is the composite of
$$\cV\pri\llra\cV\llra \Ob.
$$
In this case, the cone cycle $C\pri\sub V\pri$ is merely the pull
back of $C\sub V$ via the induced vector bundle homomorphism
$\varphi\mh V\pri\to V$.
Now let $t\pri\in \cA_\cS(\cV\pri)$ be a lifting of $s\in
\cA_\cS(\cV)$, and let $\pi\pri\mh V\pri\to W$ be the projection.
Then because $C\pri=\varphi\sta C$,
$$\pi\pri\lsta\bl\Gamma_t\cap \varphi\sta C\br\pc=
\pi\pri\lsta\bl \Gamma_{t\pri}\cap \varphi\sta
C\br\pc=\pi\lsta\bl\Gamma_s\cap C\br\pc
$$
as pseudo cycles. This shows that the pseudo-cycle representative
$\pi\bl\Gamma_s\cap C\br$ of $[W]\virt$ only depends on the section
$\xi\in\cA_\cS(\Ob)$.

This way, those sections $\xi\in\cA_\cS(\Ob)$ whose lifts intersect
transversally with the normal cone provide us pseudo-cycle
representatives of the virtual cycle $[W]\virt$. In the following,
we shall denote such representative by $D(\xi)\pc$.

In the remainder of this paper, for a scheme $W$ with obstruction
sheaf $\Ob$, we shall fix a locally free sheaf $\cO_W(V)$ that
surjects onto $\Ob$ with $C$ the virtual normal cone in $V$ of the
obstruction theory of $W$. We say that a smooth section $\xi\in
\cA_\cS(\Ob)$ is a good section if it has a lift $s\in \cA_\cS(V)$
that intersects transversally with $C$. By the previous
construction, the image $\pi\lsta(\Gamma_s\cap C)\pc$ defines a
closed pseudo-cycle $D(\xi)\pc$.

We summarize this subsection in the following Proposition.

\begin{prop}
Let the notation be as before. Then any good section
$\xi\in\cA_\cS(\Ob)$ defines a pseudo-cycle $D(\xi)$ in $W$ whose
associated homology class is the virtual cycle $[W]\virt$ in
$H_\ast(W,\ZZ)$.
\end{prop}

\subsection{Virtual cycle of Hilbert schemes of $\alpha$--points}

We shall employ smooth sections to construct cycle representatives
of the virtual cycles of $Y\hhbr{\alpha}$ for a smooth family of
quasi-projective threefolds $Y/T$ over a smooth base $T$. But before
we do that, we shall first define how cycles are pulled back by the
tautological map $Y\hhalp\to Y\hbr{\alpha}$.

We begin with defining the $\alpha$--multiplicity of the
symmetrization morphism
$$S\lalp\co  Y^\Lam\lra Y^{(\alpha)}.
$$
Let $x\in Y^\Lam$ be any element. We define the multiplicity
$m\lalp(x)$ be the number of permutations of $\Lambda$ that fix $x$
and leave $\alpha$ invariant. Namely,
$$m\lalp(x)=\#\{\sigma\in Symm(\Lam)\mid
\sigma(\alpha)=\alpha,\  \sigma(x)=x\}.
$$
It is easy to see that all elements in $(S\lalp)\upmo S\lalp(x)$
have identical $\alpha$--multiplicities; their summations
satisfies\footnote{We define $\alpha!=\alpha_1!\cdots\alpha_k!$.}
\begin{equation}\label{4.41}
\sum_{y\in(\cS^{\alpha})\upmo\cS\ualp(x)}m\lalp(y)=\alpha!.
\end{equation}
For any $z\in \hhy{\alpha}$ lies over $x\in Y^\Lam$, we define its
$\alpha$--multiplicity $m\lalp(z)=m\lalp(x)$.

Because the virtual normal cone is constructed based on the
obstruction theory of the Hilbert scheme of points $I_{Y/T}(0,n)$,
we need to work with the tautological morphism
\begin{gather*}
\varphi\lalp\co  Y\hhalp\llra I_{Y/T}(0,\alpha)
\\
\tag*{\hbox{with}}
I_{Y/T}(0,\alpha)=I_{Y/T}(0,m_1)\times_T\cdots\times_T
I_{Y/T}(0,m_k), \quad \text{where}\ m_i=|\alpha_i|.
\end{gather*}
We let $F$ be a vector bundle on $I_{Y/T}(0,\alpha)$ and let $C\sub
F$ be a cycle, which is a linear combination of subvarieties of $F$.
We let $E$ be the pull back vector bundle $\varphi\lalp\sta F$ on
$Y\hhalp$; let $\phi\lalp\mh E\to F$ and $\pi\lalp\mh E\to Y\hhalp$
be the obvious projections.

\begin{defi}\label{def4.4} For any subvariety $D\sub E$ we define the
$\alpha$--multiplicity $m\lalp(D)$ of $D$ be the
$\alpha$--multiplicity of the general point of the image
$\pi\lalp(D)\sub\hhy{\alpha}$. For any subvariety $C\sub F$ and
irreducible decomposition $\phi\lalp\upmo(C)=\cup^r_{i=1} D_i$, we
define the pull-back
$$\phi\lalp\sta C=\sum^r_{i=1} m\lalp(D_i) D_i.
$$
We define the pull back of any cycle by extension via linearity.
\end{defi}

The identity \eqref{4.41} implies that for any cycle $C$ in $F$
$$\phi_{\alpha\ast}
\phi\lalp\sta C={\alpha!}\, C.
$$
The virtual cycle of $Y\hhalp$ will be defined as the Gysin map
image of the pull back virtual normal cone given by the obstruction
theory of $I_{Y/T}(0,\alpha)$. Let $\cF\lalp$ be a locally free
sheaf on $I_{Y/T}(0,\alpha)$ that makes the obstruction sheaf of
$I_{Y/T}(0,\alpha)$ its quotient sheaf. Then the obstruction theory
of $I_{Y/T}(0,\alpha)$ provides us a cone cycle $C\lalp\in A\lsta
F\lalp$ in the vector bundle $F\lalp$ associated to $\cE\lalp$. We
let $E\lalp$ be the pull back vector bundle over $Y\hhalp$; let
$\phi\lalp\mh E\lalp\to F\lalp$ be the projection, and let
$$C\lalp=\phi\sta C\lalp
$$
be the pull back cycle in $E\lalp$ defined in \fullref{def4.4}. Since the pull back of the obstruction sheaf of
$I_{Y/T}(0,\alpha)$ is canonically isomorphic to the obstruction
sheaf of $Y\hhalp$, the sheaf $\cE\lalp=\cO(E\lalp)$ has the
obstruction sheaf $\Ob\hhalp_Y$ of $Y\hhalp$ as its quotient sheaf.
By abuse of notation, we will call the cycle $C\lalp\in A\lsta
E\lalp$ the virtual cone of the obstruction theory of $Y\hhalp$.

\begin{defi} We define the virtual fundamental class
$$[\hhy{\alpha}]\virt=0_{E\lalp}\sta \bl C\lalp\br
\in H\lsta(\hhy{\alpha},\ZZ).
$$
\end{defi}

Since the push-forward of $[\hhy{\alpha}]\virt$ under $Y\hhalp\to
I_{Y/T}(0,\alpha)$ is ${\alpha!}$ times the virtual fundamental
cycle $[I_{Y/T}(0,\alpha)]\virt$, it is independent of the choice of
$E\lalp$, thus is well-defined.

To get an explicit cycle representative, we can take a good section
$\xi\lalp$ of $E\lalp$ that intersects the cone $C\lalp$
transversally to form a pseudo-cycle
$$D(\xi\lalp)=\pi\lalp(\Gamma_{\xi\lalp}\cap C\lalp)\pc\in
PCy\lsta(\hhy{\alpha}).
$$
As was shown before, the cycle $D(\xi\lalp)$ only depend on the
image section $s\lalp\in \cA_\cS(\Ob\hhbr{\alpha}_Y)$ of $\xi\lalp$.
Hence to eliminate the dependence on $E\lalp$, we shall denote
$D(\xi\lalp)$ by $D(s\lalp)$. We have
$$[D(s\lalp)]=[\hhy{\alpha}]\virt\in H_\ast(\hhy{\alpha},\ZZ).
$$

\section{Approximation of Virtual cycles}

Because the partial equivalence $Y\hhalp\lralpbe\cong
Y\hhbe\lrbealp$ is functorial, we expect that the obstruction
sheaves, the virtual normal cones, and the cycle representatives of
the virtual cycles of $Y\hhalp$ and $Y\hhbe$ are identical over
$Y\hhalp\lralpbe\cong Y\hhbe\lrbealp$. It is the purpose of this
section to show that this is the case.

\subsection{Cones under equivalence}

Our immediate task is to compare the sheaves $\Ob\hhbr{\alpha}_Y$
with $\Ob\hhbr{\beta}_Y$ and compare their respective virtual normal
cones. Since for different $\alpha$ and $\beta$, the mentioned
partial equivalence follows from the equivalence
$Y\hhalp_{(\alpha,\alpha\wedge\beta)}\cong
Y\hhbr{\awb}_{(\awb,\alpha)}$ and
$Y\hhbe_{(\beta,\alpha\wedge\beta)}\cong
Y\hhbr{\awb}_{(\awb,\beta)}$, for our purpose we only need to
investigate the case where $\alpha>\beta$.

We first set up the notation. Let $\alpha>\beta$ be any pair. In
this section we will fix once and for all an indexing
$$\alpha=(\alpha_1,\cdots,\alpha_k)\and
\beta=(\beta_{11},\cdots,\beta_{1l_1},\cdots,\beta_{k1},\cdots\beta_{kl_k})
$$
so that $\beta_{i1}\cup\cdots\cup\beta_{il_i}=\alpha_i$. Since
$\alpha>\beta$, such indexing exists.

We let $\tilde\cZ_1,\cdots,\tilde\cZ_k$ be subschemes of
$Y\hhalp\lralpbe\times_T Y$ that are part of the universal family of
$Y\hhalp$ over $Y\hhalp\lralpbe$; we let
$\tilde\cW_{11},\cdots,\tilde\cW_{kl_k}$ be subschemes of
$Y\hhbe\lrbealp\times_T Y$ that are part of the universal family of
$Y\hhbe$ over $Y\hhbe\lrbealp$.

The partial equivalence $Y\hhalp\lralpbe\cong Y\hhbe\lrbealp$
induces a rational map from $I_{Y/T}(0,\beta)$ to
$I_{Y/T}(0,\alpha)$. Let
$$\varrho_\alpha\co  Y\hhalp\lralpbe\llra I_{Y/T}(0,\alpha)
$$
be the tautological morphism that is induced by the families
$(\cZ_1,\cdots,\cZ_k)$; we let $U\ube\lralpbe\sub I_{Y/T}(0,\alpha)$
be the image subset of this map, which is open. Because of this, we
shall endow it with the induced scheme structure (usually
non-reduced) from that of $I_{Y/T}(0,\alpha)$. For $\beta$, we have
the similarly defined
$$\varrho_\beta\co  Y\hhbe\lrbealp\llra I_{Y/T}(0,\beta)
$$
induced by the families $\cW_{ij}$. We then endow the open subset
$U\ube\lrbealp=\image(\varrho_\beta)$ with the induced scheme
structure from that of $I_{Y/T}(0,\beta)$.

More to that, for any $\tilde\eta\in U\ube\lrbealp$ over $t\in T$
that is the image of an $\eta\in Y\hhbe\lrbealp$, the associated
(indexed) zero-subschemes
$\eta_{11},\cdots,\eta_{1l_1},\cdots,\eta_{kl_k}\sub Y_t$ for each
$1\leq i\leq k$ have that the collection
$\eta_{i1},\cdots\eta_{il_i}$ is mutually disjoint. Hence we can
assign
$$\xi_i=\eta_{i1}\cup\cdots\cup \eta_{il_i},
$$
thus obtaining a zero-subscheme in $I_{Y/T}(0,|\alpha_i|)$ and the
tuple
$$(\xi_1,\cdots,\xi_k)\in I_{Y/T}(0,\alpha).
$$
It is easy to see that this correspondence defines a one--one onto
map from $U\ube\lrbealp$ to $U\ualp\lralpbe$.

The proof given in sub\fullref{sec1.3} immediately shows that

\begin{lemm}\label{iso3}
The induced map
$$\Upsilon_{\alpha\beta}\co U\ube\lrbealp\llra U\ualp\lralpbe
$$
is an \'etale morphism between two schemes; it commutes with the
maps $\varrho\lalp$, the map $\varrho\lbe$ and the partial
equivalence $Y\hhalp\lralpbe\cong Y\hhbe\lrbealp$. Further, the
obstruction sheaves $\Ob\hhalp$ of $I_{Y/T}(0,\alpha)$ are
isomorphic under pull back by $U\ube\lrbealp$:
$$\Upsilon_{\alpha\beta}\sta\Ob\hhbr{\beta}_Y
\cong \Ob_Y\hhbr{\alpha}|_{Y\ualp\lralpbe}.
$$
\end{lemm}

The virtual normal cones are also identical under this isomorphism.
We pick a locally free sheaf $\cE\lalp$ on $U\ualp\lralpbe$ that
makes $\Ob\hbr{\alpha}_Y$ its quotient sheaf (over
$U\ualp\lralpbe$). Then $\cE\lbe=\Upsilon\sta_{\alpha\beta}\cE\lalp$
is a locally free sheaf over $U\ube\lrbealp$ that makes
$\Ob_Y\hbr{\beta}$ its quotient sheaf. Then the perfect-obstruction
theory provides us the virtual normal cone $C\lalp\sub E\lalp$ and
the normal cone $C\lbe\sub E\lbe$.

\begin{lemm}\label{4.10}
Under the induced flat morphism $\phi_{\alpha\beta} \mh E\lbe\to
E\lalp$, the cycles
$$\phi_{\alpha\beta}\sta C\lalp=C\lbe.
$$
\end{lemm}

\begin{proof}
This follows from the uniqueness assertion on cones \eqref{cone1}
and the following invariance result.
\end{proof}

\begin{lemm}\label{inv5}
Let $U\sub X$ be an open subset and let $\xi\sub U$ be a zero
subscheme. Then the obstruction spaces to deforming $\xi$ in $U$ and
in $X$ are canonically isomorphic. Further, under this isomorphism
of obstruction spaces, the obstructions to deforming $\xi$ in $U$
and in $X$ are identical.
\end{lemm}

\begin{proof}
First the to obstruction spaces are traceless extension groups
$$\Ext_X^2(\cI_\xi,\cI_\xi)_0\and \Ext^2_U(\cI_\xi,\cI_\xi)_0.
$$
They are isomorphic because
$$\Ext_X^2(\cI_\xi,\cI_\xi)_0\cong H^0\bl\ext^3_{\cO_X}(\cO_\xi,\cI_\xi)\br
=H^0\bl\ext^3_{\cO_U}(\cO_\xi,\cI_\xi)\br\cong\Ext^2_U(\cI_\xi,\cI_\xi)_0.
$$
As to the obstruction theory, say using locally free resolutions of
$\cI_\xi$ and using C\v ech cohomology representative of the
obstruction classes, one checks directly that the obstruction to
deforming $\xi$ in $X$ gets mapped to the obstruction class to
deforming $\xi$ in $U$ under the canonical homomorphism
$$\Ext_X^2(\cI_\xi,\cI_\xi)_0\llra \Ext^2_U(\cI_\xi,\cI_\xi)_0.
$$
But because this arrow is an isomorphism, the two obstruction
classes must be identical.
\end{proof}

\subsection{Some further notations}

From now on, for any $\alpha$ we fix a locally free sheaf $\cE\lalp$
over $I_{Y/T}(0,\alpha)$ that makes its quotient sheaf the
obstruction sheaf $\Ob\hbr{\alpha}$ of $I_{Y/T}(0,\alpha)$; we let
$C\lalp\sub E\lalp$ be the associated virtual normal cone. Because
$Y\hbr{\alpha}$ is $I_{Y/T}(0,\alpha)$ with reduced scheme
structure, we can view $E\lalp$ as a vector bundle over
$Y\hbr{\alpha}$ and view $C\lalp$ as a cone cycle in $A\lsta
E\lalp$.

For the Hilbert scheme of $\alpha$--points $Y\hhalp$, its obstruction
sheaf $\Ob\hhalp_Y$ is canonically isomorphic to the pull back sheaf
$\rho\lalp\sta\Ob\halp\sub I_{Y/T}(0,\alpha)$ under the tautological
morphism
$$\rho\lalp\mh Y\hhalp\to Y\halp.
$$
Thus by taking $\tilde\cE\lalp=\rho\lalp\sta \cE\lalp$, the
obstruction sheaf $\Ob_Y\hhalp$ naturally becomes a quotient sheaf
of $\tilde\cE\lalp$.

We let $\tilde E\lalp\to E\lalp$ be the tautological projection,
viewed as a stack flat morphism
$$\tilde E\lalp= E\lalp\times_{Y^{(\alpha)}}Y\ualp\llra E^\Lam.
$$
The normal cone $\tilde C\lalp\sub\tilde E\lalp$ is then defined to
be the stack flat pull back of $C\lalp$ as specified in \fullref{def4.4}.

For stratifications of $\tilde C\lalp$ and $Y\hhalp$, we shall take
the standard pair of stratifications that respects the morphism
$\tilde C\lalp\to Y\hhalp$ and the loci of non-locally freeness of
the sheaf $\Ob_Y\hhalp$.

\begin{defi}
We say that a smooth section $s\lalp\sub \cA(\Ob_Y\hhalp)$
intersects transversally with its normal cone if one (thus all) of
its lifts $\xi\lalp\in \cA(\tilde\cE\lalp)$ intersects transversally
with the cone $\tilde C\lalp$.
\end{defi}

Following the discussion in sub\fullref{rep4.4}, we can find
sections of $\cA(\Ob_Y\hhalp)$ that intersect the normal cone
transversally. For such $s\lalp$, we shall denote by $D(s\lalp)\sub
Y\hhalp$ the pseudo-cycle that is the image in $Y\hhalp$ under
$\tilde E\lalp\to Y\hhalp$ of intersecting a lift of $s\lalp$ with
$\tilde C\lalp\sub\tilde E\lalp$.

Now let $s\lalp\in \cA(\Ob_Y\hhalp)$ be a section that intersects
transversally with the normal cone. Because
$\Ob_Y\hhalp|_{Y\hhalp\lralpbe}\cong \Ob_Y\hhbe|_{Y\hhbe\lrbealp}$
canonically, we can view
$s_{\alpha|\beta}=s\lalp|_{Y\hhalp\lralpbe}$ as a section of
$\Ob_Y\hhbe$ over $Y\hhbe\lrbealp$. Because of \fullref{iso3} and
\ref{4.10}, $s_{\alpha|\beta}$ is a smooth section of $\Ob_Y\hhbe$
and intersects transversally with the normal cone of $Y\hhbe$.

\subsection{Compatible cycle representatives}

To compare the cycles $[Y\hhbr{\alpha}]\vir$, in this section we
shall carefully pick smooth sections $s\lalp$ so that for any pair
$\alpha$ and $\beta$ the cycle representatives $D(s\lalp)$ are
$D(s\lbe)$ are identical over most part of the intersection
$\hhy{\alpha}_{(\alpha,\beta)}\cong Y\hhbe\lrbealp$.

To this end, we form the strict $\alpha$--diagonal
$$\Delta\lalp=\{x\in Y^\Lambda\mid a\sim_\alpha b \Rightarrow
x_a=x_b\};
$$
they are closed; for different $\alpha$ and $\beta$,
$\Delta\lalp\cap\Delta\lbe=\Delta_{\alpha\vee\beta}$\footnote{$\alpha\vee\beta$
is the smallest element among all that are larger than or equal to
both $\alpha$ and $\beta$.}. Next we fix a sufficiently small $c>0$
and pick a function $\varepsilon\mh \cP\llam\to (0,c)$ whose values
on any ordered pair $\alpha>\beta$ obey
$\varepsilon(\alpha)>R\varepsilon(\beta)$ for a sufficiently large
$R$. After fixing a Riemannian metric on $Y$, we then form the
$\varepsilon$--neighborhoods of $\Delta\lalp\sub Y^\Lambda$:
$$\Delta_{\alpha,\varepsilon}=\{x\in Y^\Lambda\mid \text{dist}(x,\Delta\lalp)<
\varepsilon(\alpha)\}.
$$
For any $\beta\leq\alpha$, we form
$$\Delta\ualp_{\beta,\varepsilon}=\bigcup_{\alpha\geq\gamma\geq\beta}
\Delta_{\gamma,\varepsilon} \and
Q\ualp_{\beta,\varepsilon}=\Delta_{\beta,\varepsilon}-\bigcup_{\alpha\geq\gamma>\beta}
\Delta_{\gamma,\varepsilon}\ualp.
$$
Note that $Q\ualp_{\beta,\varepsilon}$ are closed subsets of
$\Delta_{\beta,\varepsilon}$.

We have the following intersection property of these sets.

\begin{lemm}\label{int}
For any pair $\beta_1,\beta_2\leq\alpha$ satisfying
$\Delta_{\beta_1,\varepsilon}\cap
Q\ualp_{\beta_2,\varepsilon}\ne\emptyset$, necessarily $\beta_2\geq
\beta_1$.
\end{lemm}

\begin{proof}
Because $c$ is sufficiently small, whenever
$\Delta_{\mu,\varepsilon}\cap\Delta_{\nu,\varepsilon}\ne \emptyset$,
necessarily $\Delta_{\mu}\cap\Delta_{\nu}\ne \emptyset$. Then
because $\Delta_{\mu}\cap\Delta_{\nu}=\Delta_{\mu\vee\nu}$, because
$\Delta_\mu$ intersects $\Delta_{\nu}$ perpendicularly, and because
$\varepsilon(\mu\vee\nu)>{\frac{R}{2}}\varepsilon(\mu)+{\frac{R}{2}}\varepsilon(\nu)$,
$$\Delta_{\mu,\varepsilon}\cap \Delta_{\nu,\varepsilon}\sub
\Delta_{\mu\vee\nu,\varepsilon}.
$$
Now suppose $\beta_2\ngeq\beta_1$, then
$\beta_1\vee\beta_2>\beta_2$; therefore we have
$\Delta_{\beta_1,\varepsilon}\cap\Delta_{\beta_2,\varepsilon}\sub\Delta_{\beta_1\vee\beta_2,\varepsilon}$
and $Q\ualp_{\beta_2,\varepsilon}\sub
\Delta\ualp_{\beta_2,\varepsilon}-\Delta\ualp_{\beta_1\vee\beta_2,\varepsilon}$.
Combined, we have
$$\emptyset\ne
\Delta_{\beta_1,\varepsilon}\cap Q\ualp_{\beta_2,\varepsilon}\sub
\Delta_{\beta_1\vee\beta_2,\varepsilon}\cap\bl
\Delta_{\beta_2,\varepsilon}-
\Delta_{\beta_1\vee\beta_2,\varepsilon}\br=\emptyset.
$$
This proves $\beta_2\geq\beta_1$.
\end{proof}

An immediate corollary of this is that
$\Delta\ualp_{\beta,\varepsilon}=\coprod_{\alpha\geq\gamma\geq\beta}
Q\ualp_{\gamma,\varepsilon}$ forms a partition (a disjoint union) of
$\Delta\ualp_{\beta,\varepsilon}$. By choosing $\beta=0_\Lambda$, it
also shows that $\{Q\ualp_{\beta,\varepsilon}\mid \beta\leq\alpha\}$
forms a partition of $Y^\Lambda$.

Moving to $\hhy{\alpha}$, we form
$$\cN^{\alpha}_{\beta,\varepsilon}=\rho\lalp\upmo\bl
N_{\beta,\varepsilon}\ualp\br\and
\cQ^{\alpha}_{\beta,\varepsilon}\defeq \rho\lalp\upmo\bl
Q\ualp_{\beta,\varepsilon}\br \ \sub \hhy\alpha,
$$
in which we continue to denote by $\rho\lalp\mh Y\hhbr{\alpha}\lra
Y^\Lam$
the projection. The collection $\{\cQ^\alpha_{\beta,\varepsilon}\mid
\beta\leq\alpha\}$ forms a partition of $\hhy\alpha$.

\eject

\begin{lemm}\label{req} For sufficiently small $c$, we can find a collection
of sections $s\lalp\in \cA_\cS(\Ob\hhbr{\alpha})$ that satisfy the
properties
\begin{itemize}
\item[\rm(i)] {\sl each $s\lalp$ intersects transversally
with the normal cone of $\Ob_Y\hhalp$;}

\item[\rm(ii)] {\sl for any $\beta<\alpha$, the sections $s\lalp$ and $s\lbe$
coincide over $\cQ^{\alpha}_{\beta,\varepsilon}$}.
\end{itemize}
\end{lemm}

The requirement (ii) is understood as follows. To each
$\beta<\alpha$, because $\cQ^{\alpha}_{\beta,\varepsilon}$ is
disjoint from $\rho\lalp\upmo(\Delta_\gamma)$ for all
$\alpha\geq\gamma>\beta$, it lies inside $\hhy{\alpha}\lralpbe$.
Thus restricting to $\cQ^{\alpha}_{\beta,\varepsilon}$
both $s\lalp$ and $s\lbe$ are sections of the same sheaf, and hence
can be said to equal.

\begin{proof}
We prove the lemma by induction. Because the space $\hhy{\alpha}$ is
a disjoint union of $\cQ^\alpha_{\beta,\varepsilon}$, we will
construct $s\lalp$ by specifying its values along each of the above
subsets according to (ii) and then showing that the resulting
section can be extended to satisfy (i).

We now construct the section $s\lalp$ by induction. Suppose we have
already constructed $s\lbe$ for all $\beta<\alpha$ that satisfy the
properties (i)--(ii). Along the partition
$Y\hhbr{\alpha}=\coprod_{\beta\leq\alpha}\cQ^{\alpha}_{\beta,\varepsilon}$,
we shall follow the rule (ii) to define
$$s_{\alpha|\beta}=s\lbe|_{\cQ\ualp\lrbealp}
\in\Gamma(\cQ^{\alpha}_{\beta,\varepsilon},\cA_\cS(\Ob_Y\hhbr{\alpha})
$$
be the restriction to $\cQ^{\alpha}_{\beta,\varepsilon}$ of $s\lbe$.
Inductively, this will define
$$s\lalp\in\Gamma\bl\hhy{\alpha}-\cN^{\alpha}_{\alpha,\varepsilon},
\cA_\cS(\Ob_Y\hhbr{\alpha})\br
$$
after checking that the collection $s_{\alpha|\beta}$ forms a smooth
section over $\hhy{\alpha}-\cN^{\alpha}_{\alpha,\varepsilon}$.

We now prove that it is so. First, because
$\hhy{\alpha}-\cN^{\alpha}_{\alpha,\varepsilon}$ is a disjoint union
of $\{\cQ^{\alpha}_{\beta,\varepsilon}\mid \beta<\alpha\}$, each
$z\in \hhy{\alpha}-\cN^{\alpha}_{\alpha,\varepsilon}$ must lie in a
$\cQ^{\alpha}_{\beta,\varepsilon}$ for a unique $\beta<\alpha$. In
case $z$ is an interior point of $\cQ^{\alpha}_{\beta,\varepsilon}$,
then $s\lalp$ coincide with $s\lbe$ near $z$; by induction
hypothesis, $s\lalp$ satisfies the requirement (i) and (ii) near
$z$. In case $z$ is not an interior point of
$\cQ^{\alpha}_{\beta,\varepsilon}$, since
$\cN^{\alpha}_{\gamma,\varepsilon}$ is open, by \fullref{int} this
is possible only when $z$ lies in the closure $cl\bl
\cQ^{\alpha}_{\beta\pri,\varepsilon}\br$ for some $\beta\pri>\beta$.
There are two possibilities: one is when $\beta\pri=\alpha$, in
which case nothing to prove. The other is when $\beta\pri<\alpha$.
Since $s\lalp$ is defined via $s\lbe$ on
$\cQ^{\alpha}_{\beta,\varepsilon}$ and via $s_{\beta\pri}$ on
$\cQ^{\alpha}_{\beta\pri,\varepsilon}$, to check the continuity of
$s\lalp$, we need to compare the germ of $s\lbe$ and of
$s_{\beta\pri}$ near $z$. In this case,
$\cQ^{\alpha}_{\beta\pri,\varepsilon}\sub Y\hhalp\lralpbe$, thus can
be considered lies in $Y\hhbr{\beta\pri}$; as to
$\cQ\ualp_{\beta,\varepsilon}$,
$\cQ^{\alpha}_{\beta,\varepsilon}\cap
Y\hhbr{\beta\pri}=\cQ^{\beta\pri}_{\beta,\varepsilon}$. Hence by
induction hypothesis, $s_{\beta\pri}$ is an extension of $s\lbe$,
Thus $s\lalp$ is well-defined near $z$, and thus is smooth near $z$.

Once we know that $s\lalp$ is smooth over
$\hhy\alpha-\cN^{\alpha}_{\alpha,\varepsilon}$, which is closed in
$Y\hhbr{\alpha}$, we can extend it to a good smooth section of
$\Ob_Y\hhbr{\alpha}$ over $Y\hhbr{\alpha}$, which completes the
proof of the Lemma.
\end{proof}

From now on, we fix such a collection
$\{s\lalp\}_{\alpha\in\cP\llam}$ and form their associated
pseudo-cycle representatives $D(s\lalp)$.

\subsection{Approximation of the virtual cycles}

In this subsection, we shall investigate how pseudo-cycles
$D(s\lalp)$ are related by looking at their images in $Y^\Lam$.

We let $d$ be the real dimension of $T$; let $D(s\lalp)$ be the
$d$--dimensional pseudo-cycles constructed relative to a set $\tilde
\Theta\lalp\sub Y\hhalp$ of dimension $\leq d-2$. We take
$\Theta\sub Y^\Lam$ be the union of the images under the projection
$\rho\lalp\mh Y\hhalp\to Y^\Lam$ of all $\tilde \Theta\lalp$:
$\Theta=\cup\lalp\rho\lalp(\tilde\Theta\lalp)$. Then to a
sufficiently small $\eps>0$ independent of $c$, we pick $\tilde
f\lalp\mh \Sigma\lalp\to Y\hhalp$ a representative of $D(s\lalp)$ as
$d$--dimensional pseudo-cycle relative to
$(\rho\lalp\upmo(\Theta),\eps)$. The composition
$f\lalp=\rho\lalp\circ\tilde f\lalp$ defines a
$(\Theta,\eps)$--relative pseudo-cycle, which we denote by
$[f\lalp]\pc$.

We now relate different pseudo-cycles $[f\lalp]\pc$ by forming
inductively
\begin{equation}
\label{def-2} \delta\lalp=[f\lalp]\pc-\sum_{\beta<\alpha}\delta\lbe\in
Cyc_d(Y^\Lam)_\Theta.
\end{equation}
We have the following vanishing result:

\begin{lemm} \label{lem5.6}
The pseudo-chain $\delta\lalp\cap(Y^\Lam-\Delta_{\alpha,2c})$ is
equivalent to $0$.
\end{lemm}

\begin{proof}
What we need to show is that to any $z\in
Y^\Lam-(\Delta_{\alpha,2c}\cup\Theta_\epsilon)$ we can find a
sufficiently small ball $B_r(z)$ centered at $z$ so that
$\delta\lalp\cap B_r(z)$ as a pseudo-cycle is equivalent to zero.

We now prove this by induction. For $\alpha=0\llam$, there is
nothing to prove. Now let $\alpha\in\cP\llam$ be any element so that
this holds true for all $\beta<\alpha$. Let $x\in Y^\Lam$ be any
point away from $\Delta_{\alpha,2c}\cup\Theta_\epsilon$ and let
$\beta<\alpha$ be so that
\begin{equation}\label{111}
x\in
Q\ualp_{\beta,\varepsilon}.
\end{equation}
We first claim that for sufficiently small $r$, $\delta_\gamma\cap
B_r(x)\sim 0$ for all $\beta\not\geq\gamma<\alpha$. Suppose not,
then by induction hypothesis, $x\in \Delta_{\gamma,2c}$. Since $c$
is sufficiently small and
$\varepsilon(\beta_1)>R\varepsilon(\beta_2)$ for any
$\beta_1>\beta_2$, we have
$\Delta_{\beta,\varepsilon}\cap\Delta_{\gamma,\varepsilon}\sub
\Delta_{\beta\vee\gamma,\varepsilon}$, and hence $\rho(z)\in
\Delta_{\beta\vee\gamma,\varepsilon}$. Because we have assumed that
$\beta\not\geq\gamma$, we must have $\beta\vee\gamma> \beta$, which
implies that $\Delta_{\beta\vee\gamma,\varepsilon}\cap
Q\ualp_{\beta,\varepsilon}=\emptyset$, contradicting to $x\in
Q\ualp_{\beta,\varepsilon}$ and in
$\Delta_{\beta\wedge\gamma,\varepsilon}$.

Because of this, when we intersects $\delta\lalp$ with $B_r(x)$, for
sufficiently small $r$ \eqref{def-2} tells us that
$$\delta\lalp\cap B_r(x)=
[f\lalp]\pc\cap B_r(x)-\sum_{\gamma\leq\beta}\delta_\gamma\cap
B_r(x).
$$
On the other hand by definition
$$\delta\lbe\cap B_r(x)=[f\lbe]\pc\cap B_r(x)-
\sum_{\gamma<\beta}\delta_\gamma\cap B_r(x).
$$
Hence if we can show that $[f\lalp]\pc\cap B_r(x)=[f\lbe]\pc\cap
B_r(x)$, the above two identities will force $\delta\lalp\cap
B_r(x)=0$, exactly what we intend to prove.

For this we argue as follows: let $z\in
D(s\lalp)\cap\rho\lalp\upmo(B_r(x))$. Because of \eqref{111}, $x\in
\Delta_{\beta,\varepsilon}$ and
$x\not\in\cup_{\beta<\gamma\leq\alpha}\Delta_\gamma$; thus $z$ must
lie in $\hhy\alpha\lralpbe$ and thus also in $\hhy\beta\lrbealp$.
Hence (ii) of lemma \eqref{req} implies that
$$D(s\lalp)\cap\rho\lalp\upmo(B_r(x))=D(s\lbe)\cap\rho\lbe\upmo(B_r(x)).
$$
For the same reason, the above identity holds in case $x\in
D(s\lbe)$. This proves the lemma.
\end{proof}

\subsection{Truncated discrepancy cycles}

Because $\delta\lalp$ is equivalent to zero away from
$\Delta_{\alpha,2c}$, we can truncate it by intersecting it with the
closed subset $\overline\Delta_{\alpha,e}$ with a general $2c<e<3c$:
$\delta\lalp^{e}=\delta\lalp\cap\overline\Delta_{\alpha,e}$. It can
also be defined inductively by
$$\delta\lalp^{e}\sim [f\lalp]\pc\cap\overline\Delta_{\alpha,e}
-\sum_{\beta<\alpha}\delta\lbe^{e},
$$
where the intersection $[f\lalp]\pc\cap\overline\Delta_{\alpha,e}$
can be replaced by any pseudo-chain $f\lalp|_{\Sigma\lalp^{e}}$ with
$\Sigma\lalp^{e}=f\lalp\upmo(\overline\Delta_{\alpha,e})-O$ for an
open $O$ such that $\Sigma\lalp^{e}$ has smooth boundary and
$f\lalp(O)\sub\Theta_\varepsilon$. By Sards theorem, for general $e$
the set $f\lalp\upmo(\overline\Delta_{\alpha,e})$ has smooth
boundary away from $f\lalp\upmo(\Theta_\varepsilon)$. Thus we can
choose $e$ that works for all $\alpha$.

\begin{lemm}
With such choice of $e$ and pseudo-chain representative
$[f\lalp]\pc\cap\overline\Delta_{\alpha,e}$, the inductively defined
pseudo chain $\delta\lalp^{e}$ is a pseudo-cycle relative to
$(\Theta,\epsilon)$.
\end{lemm}

\begin{proof}
This follows directly from \fullref{lem5.6} and the definition of
pseudo-cycle.
\end{proof}

In case $\alpha$ has more than one equivalence class, the cycle
$\delta\lalp$ is equivalent to the product of $\delta_{\alpha_i}$.
More precise, by viewing each $\alpha_i$ as a set, we can form the
spaces $Y^{\alpha_i}$ and the pseudo-cycles $\delta_{\beta_i}$ for
$\beta_i\in\cP_{\alpha_i}$. We continue to denote by ${\alpha_i}$
the top partition in $\cP_{\alpha_i}$.

\begin{lemm}\label{prod1}
We can choose representatives $\delta_{\alpha_i}$ and
$\delta_\alpha$ so that under the identity $Y^\Lam=\prod_T^k
Y^{\alpha_i}$ as pseudo-cycles $\delta\lalp\cong \prod_{i=1}^k
\delta_{{\alpha_i}}$.
\end{lemm}

\begin{proof}
For each $\alpha_i$, we form the cycle representatives
$D(\gamma)\sub Y\hhbr{\gamma}$ and the discrepancy cycles
$\delta_{\gamma}$ by picking a collection of sections
$\{s_{\gamma}\mid\gamma\in\cP_{\alpha_i}\}$ provided by the previous
lemma. For each $\beta\leq\alpha$ in $\cP\llam$, we write $\beta$ as
$(\beta_{11},\cdots\beta_{kl_k})$ and form
$\beta_i=(\beta_{i1},\cdots,\beta_{il_i})$, each is a partition in
$\cP_{\alpha_i}$. Then
$$\hhy\beta=\hhy{\beta_1}\times\cdots\times\hhy{\beta_k}
\and \Ob_Y\hhbr\beta=\pi_1\sta\Ob_Y\hhbr{\beta_1}
\oplus\cdots\oplus\pi_k\sta\Ob_Y\hhbr{\beta_k}
$$
with $\pi_i$ the $i$-th projection of $Y\hhalp$ to
$Y\hhbr{\alpha_i}$. Further, the sections
$s_{\beta_1},\cdots,s_{\beta_k}$ provides us a section
$$s\lbe=\pi_1\sta s_{\beta_1}
\oplus\cdots\oplus\pi_k\sta s_{\beta_k}
$$
with associated pseudo-cycle
\begin{equation}\label{4.1}
D(s\lbe)=D(s_{\beta_1})\times\cdots\times D(s_{\beta_k}).
\end{equation}
We claim that using such decomposition, the discrepancy
pseudo-cycles
\begin{equation}\label{4.2}
\delta\lbe\cong\delta_{\beta_1}\times\cdots\times\delta_{\beta_k}.
\end{equation}
Indeed, when $\beta=0\llam$, then the identity reduces to
$$D(s\lbe)=\prod_{a\in\Lambda} D(s_{1_{\{a\}}}),
$$
which is \eqref{4.1}. Now suppose the identity holds for all
$\gamma<\beta$. Then
$$\prod_{i=1}^k D(s_{\beta_i})=\prod_{i=1}^k
\Bigl(\sum_{\beta_{ij}\leq\beta_i}\delta_{\beta_{ij}}\Bigr)
=\sum_{(\gamma_1,\cdots,\gamma_k)\in\cP_{\beta_1}\times\cdots
\times\cP_{\beta_k}}\delta_{\gamma_1}\times\cdots\times\delta_{\gamma_k}
$$
$$\qquad
=\sum_{\gamma<\beta}\delta_\gamma+
\delta_{1_{\beta_1}}\times\cdots\times\delta_{1_{\beta_k}}.
\qquad\qquad\qquad\qquad\quad\
$$
The desired identity \eqref{4.2} then follows from
$D(s\lbe)=\sum_{\gamma<\beta}\delta_\gamma+\delta_\beta$ and the
identity \eqref{4.1}. This proves the Lemma.
\end{proof}

By choosing a general $e$ as before, the truncated discrepancy cycle
also satisfies the product formula \eqref{4.2} with $\delta_{\cdot}$
replaced by its truncated version $\delta_{\cdot}^e$.

\section{Proof of the Main theorem}

The proof of the theorem now is fairly straightforward. We first
apply the previous construction to $Y=X$ and to the set $\Lam=[n]$
of integers from $1$ to $n$ to relate the degree of the virtual
cycle
$$\deg[I_X(0,n)]\virt=\frac{1}{n!}\deg [X\hhn]\virt.
$$
(Here we follow the convention $X\hhbr{n}=X\hhbr{\Lam}$.)

The right hand side can be expressed as the degree of an explicit
zero-cycle. For this we pick a smooth section $s_{[n]}$ of the
obstruction sheaf $\Ob_X\hhn$ of $X\hhn$ that intersects
transversally with the normal cone. Because the virtual dimension of
$I_X(0,n)$ is zero, the resulting pseudo-cycle $D(s_{[n]})$ in
$X\hhn$ is a zero-cycle, satisfying
$$\deg [X\hhn]\virt=\deg [D(s_{[n]})].
$$
To proceed, we shall relate it to the degrees of the discrepancy
cycles $\delta\lalp$. For a sufficiently small $\eps$ and $c$ and
for all $\alpha\in\cP_\Lam$ we choose sections $s\lalp$ of
$\Ob\hhalp_X$ according to \fullref{req}. We let $D(s\lalp)$ in
$X\hhalp$ be the associated zero-cycle derived by intersecting the
normal cone by $s\lalp$. However, because each $\delta\lalp$ is a
zero-cycle, \fullref{lem5.6} shows that it is entirely contained
in the $\eps$--tubular neighborhood of $X_\Del^n\sub X^n$. Further,
their degrees
$$\deg [D(s_{[n]})]=\sum_{\alpha\in\cP_\Lam}\deg \delta\lalp.
$$
Thus to prove that $\deg[I_X(0,n)]\virt$ is expressible in terms of
a universal polynomial in Chern numbers of $X$ we only need to show
that the same hold for all $\deg\delta\lalp$, which then follows
from that of $\deg\delta_{[n]}$ by the product formula in \fullref{prod1}. 
This way, we are reduced to show that $\deg \delta_{[n]}$ is
expressible by a universal expression in Chern numbers of $X$.

Our next step is to use the construction of the diffeomorphism of a
neighborhood of $0_{TX}\sub TX$ with a neighborhood of the diagonal
$\Delta(X)\sub X\times X$ to transfer the cycle $\delta\lon$ to the
tangent bundle $TX$.

We let $\cU\sub X\times X$ be the tubular neighborhood of
$\Del(X)\sub X\times X$, let $\cV\sub TX$ be the tubular
neighborhood of $0_{TX}\sub TX$ and let $\psi\mh \cV\to\cU$ be the
smooth isomorphism provided by \fullref{lem2.1}. The map $\psi$
induces a smooth isomorphism
\begin{equation}\label{14.1}
\Psi_\alpha\co 
\cvalp_0\llra U\hhbr{\alpha}\sub X\hhalp
\end{equation}
onto an open neighborhood of $X_\Del\hhalp\sub X\hhalp$ (see \fullref{le2.6}). Because $U\hhalp$ is an open neighborhood of
$X_\Del\hhalp$, by choosing $c$ sufficiently small, \fullref{lem5.6} tells us that the cycle $\delta\lon$ lies entirely in
$U\hhn$. Thus the degrees of
$$\delta_{[n],\cV}=\Psi_{[n]}\sta(\delta\lon)
$$
is identical to that of $\delta\lon$.

The cycle $\delta_{[n],\cV}$ can also be constructed using the
universal family over the Grassmannian of quotients $\CC^N\to
\CC^3$. Let $Q\to Gr$ be the total space of the universal quotient
bundle and let $Q\hhalp$ (resp. $Q\hhalp_0$) be the relative Hilbert
scheme of $\alpha$--points (resp. centered $\alpha$--points) of
$Q/Gr$. To each $Q\hhalp_0$ we form its obstruction sheaf
$\Ob\hhalp_{Q}$; we pick a locally free sheaf $\cE\lalp$ making
$\Ob\hhalp_Q$ its quotient sheaf; we let $E\lalp$ be the associated
vector bundle and let $D\lalp\sub E\lalp$ be the associated normal
cone.

We then pick a smooth section $t\lalp$ of $E\lalp$ for all $\alpha$
so that the collection $\{t\lalp\}_{\alpha\in\cP_\Lam}$ satisfies
the conclusion of \fullref{lem2.1}. The sections $t\lalp$ provide
us the pseudo-cycle $D(t\lalp)$ in $Q\hhalp_0$ relative to a subset
$\tilde\Theta\lalp$ of dimension at most $2\dim Gr-8$. According to
\eqref{def-2}, the images under the projections $\eta\lalp\mh
Q_0\hhalp\to Q_0^n$ of the pseudo-cycles $D(t\lalp)$ form the
discrepancy pseudo-cycle $\delta_{[n],Q}\sub Q^n_0$ relative to
$(\Theta,\eps)$, the set $\Theta$ which is the union of all
$\eta\lalp(\tilde\Theta\lalp)$ and for $\eps$ which is sufficiently
small.

The cycle $\delta_{[n],Q}$ defines a codimension six homology class
in $Gr$. Indeed, by \fullref{lem5.6}, the pseudo-cycle is zero
outside the $2c$--tubular neighborhood in $Q_0^n$ of the top diagonal
$Q^n_\Del\cap Q_0^n$. Because the top diagonal $Q^n_\Del$ consists
of points $(\xi,\cdots,\xi)$, they are centered only if $\xi=0$.
Thus the top diagonal is the zero section of $Q_0^n/Gr$; thus is
compact. Therefore, $\delta_{[n],Q}$ is a compact pseudo-cycle, thus
defines a homology class
$$[\delta_{[n],Q}\in H_{\dim_\RR Gr-6}(Q_0^n,\ZZ)\cong
H_{\dim_\RR Gr-6}(Gr,\ZZ).
$$
We shall relate the cycle $\delta_{[n],Q}$ with $\delta_{[n]}$ in
$X^n$ by picking a smooth map $g\mh X\to Gr$ so that as smooth
vector bundle $g\sta Q=TX$. Without loss of generality, we can
choose $N$ to be sufficiently large and choose $g$ to be an
embedding. Obviously, $g$ induces a smooth isomorphism (compare to
\eqref{8.2})
\begin{equation}\label{6.1}
g_n\co  Q_0^n\times_{Gr}X\cong (TX)_0^n.
\end{equation}
We have the following compatibility result.

\begin{lemm}\label{6.61}
We can choose $g$ and sections $s\lalp$ so that $Q_0^n\times_{Gr}X$
intersects transversally with $\delta_{[n],Q}$ and the intersection
$$g_n\sta(\delta_{[n],Q})\defeq \delta_{[n],Q}\cap \bl
Q_0^n\times_{Gr}X\br,
$$
considered as a cycle in $(TX)_0^n$ via the isomorphism \eqref{6.1},
is identical to the cycle $\delta_{[n],\cV}$ via the inclusion
$\cV_0^n\sub (TX)_0^n$.
\end{lemm}

Once the lemma is proved, then we can use the constructed sections
$s\lalp$ to form the cycle $\delta\lon$ in $X^n$ to conclude
$$\deg\delta\lon
=\deg\delta_{[n],\cV}=\deg g_n\sta(\delta_{[n],Q}) =\deg\bl
[\delta_{[n],Q}]^{P.D.}\cap g\lsta([X])\br,
$$
where $[\delta_{[n],Q}]^{P.D.}$ is the Poincare dual of the homology
class $[\delta_{[n],Q}]\in H_6(Gr,\ZZ)$. Because $H\sta(Gr,\ZZ)$ is
generated by $u_{i}\in H^{2i}(Gr,\ZZ)$,
$$[\delta_{[n],Q}]^{P.D.}=P_n(u_1,u_2,u_3)\in H^6(Gr,\ZZ)
$$
expressible in a polynomial in $u_i$. On the other hand, using the
embedding $Gr(N,3)\sub Gr(N+1,3)$ and the proof that will follow, we
see immediately that this polynomial $P_n$ is independent of the
choice of $N$; thus is universal in $n$.

Finally, because $g_n\sta(u_i)=c_i(X)$,
$$[\delta_{[n],Q}]^{P.D.}\cap
g\lsta([X])=p_n(c_1(X),c_2(X),c_3(X))\cap [X]
$$
is a universal expression in the Chern numbers of $X$. This proves
the main theorem.

We shall divide the proof of \fullref{6.61} into two steps: one is
to compare the intersection $g\sta_n(\delta_{[n],Q})$ with a
similarly constructed cycle $\delta_{[n],TX}$ on $(TX)^n_0$; the
other is to compare the later with the pull back of $\delta\lon$
using the map $\cV_0^n\to X^n$. Since the proofs of both are
similar, we shall provide the details of the first while indicating
the necessary changes required for the second.

We begin with a quick account of the moduli of $\At$. We let $V=\At$
and let $V\hhalp_0\to V^n_0$ be the tautological map from the
Hilbert scheme of centered $\alpha$--points to $V_0^n$ that was
constructed in section 1 with $Y$ replaced by $V$. Because $GL(3)$
acts on $V$, it acts on $V\hhalp$ and $V^n$. Further these actions
leave the projection $V\hhalp\to V^n$ and the averaging map $V^n\to
V$ invariant, thus the space $V\hhalp_0= V\hhalp\times_V 0$ is
$GL(3)$--invariant. Not only that, its obstruction sheaf
$\Ob\hhalp_V$ and its obstruction theory are all naturally
$GL(3)$--linearized.

Because $GL(3)$ is reductive, we can find a $GL(3)$--linearized
locally free sheaf $\cF\lalp$ on $V\hhalp_0$ that makes
$\Ob\hhalp_V$ its $GL(3)$--equivariant quotient sheaf
\begin{equation}\label{12.1}
\cF\lalp\llra \Ob\hhalp_V.
\end{equation}
We let $F\lalp$ be its associated vector bundle and let $D\lalp\sub
F\lalp$ be the associated normal cone. Then each $D\lalpi$ in the
irreducible decomposition
$$D\lalp=\sum_i m_{\alpha,i} D_{\alpha,i}
$$
is a $GL(3)$--invariant cone-like subvariety of $F\lalp$.

We next fix a standard stratification of $\cup_i D\lalpi\to
V\hhalp_0$ subordinating to the loci of the non-locally freeness of
the sheaf $\Ob\hhalp_V$. To each $D\lalpi$, we let $S\lalpi\sub
D\lalpi$ be its open stratum and let $T\lalpi\sub V_0\hhalp$ be the
image of $S\lalpi$, which is a stratum of $V\hhalp_0$. Because
$D\lalpi\to F\lalp\to V\hhalp_0$ are $GL(3)$--equivariant, every
stratum, including $S\lalpi$, are $GL(3)$--invariant.

Our next step is to use \eqref{12.1} to build locally free sheaves
on $Q\hhalp_0$ and on $(TX)\hhalp_0$ making their respective
obstruction sheaves their quotient. For $Q/Gr$, we first cover $Gr$
by open $U_a\sub Gr$ with vector bundle isomorphisms
\begin{equation}
f_a\co  Q\times_{Gr}U_a\twomapright{\cong}V\times U_a.
\end{equation}
Then using the induced isomorphism
$$Q\hhalp_0\times_{Gr} U_a\cong V\hhalp_0\times U_a
$$
and $p_V$ the first projection of the product on the right hand
side, we can form the induced
\begin{equation}\label{quot1}
\cE_{\alpha,a}=p_V\sta\cF\lalp\llra p_V\Ob\hhalp_V\cong
\Ob\hhalp_Q|_{Q\hhalp_0\times_{Gr} U_a}.
\end{equation}
Over $U_{ab}=U_a\cap U_b$, the isomorphisms
$$f_{ba}=f_b\circ f_a\upmo\co  V\times U_{ab}\llra V\times U_{ab}
$$
induce transition isomorphisms
$$f_{ba}\ualp\co \cE_{\alpha,a}|_{Q\hhalp_0\times_{Gr}U_{ab}}\llra
\cE_{\alpha,b}|_{Q\hhalp_0\times_{Gr} U_{ab}}.
$$
Since $f_{ab}$ satisfy the cocycle condition, $f\ualp_{ab}$ also
satisfy the cocycle condition. Hence $f_{ab}\ualp$ glue to form a
locally free sheaf $\cE\lalp$ on $Q\hhalp_0$. Obviously, the
quotient homomorphisms \eqref{quot1} glue to form a quotient
homomorphism
$$\cE\lalp\lra \Ob\hhalp_Q.
$$
The normal cone of $\Ob\hhalp_Q$ takes a simple form in this
setting. Let $E\lalp$ be the associated vector bundle of $\cE\lalp$.
The cycles
$$D\lalp\times U_a=\sum_i m\lalpi D\lalpi\times U_a\in C\lsta( F\lalp\times
U_a)=C\lsta( E\lalp|_{Q\hhalp_0\times_{Gr}U_a})
$$
glue together to form a cycle that is the normal cone $C\lalp\sub
E\lalp$:
$$C\lalp=\sum_i m\lalpi C\lalpi
$$
with $C\lalpi$ the glued subvariety from $D\lalpi\times U_a$.

Similarly, for the vector bundle $TX/X$ we can carry over the same
procedure to form a locally free sheaf $\tilde \cE\lalp$ over
$(TX)\hhalp_0$ making the obstruction sheaf $\Ob\hhalp_{TX}$ its
quotient sheaf; the normal cone in the associate vector bundle
$\tilde C\lalp\sub\tilde E\lalp$ is also the similarly induced cycle
by $D\lalp\sub F\lalp$.

Now we look at the smooth map $g\mh X\to Gr$, which we assume to be
an embedding, and the smooth isomorphism
\begin{equation}\label{iso4}
g_0\co TX\cong Q\times_{Gr}X\sub Q.
\end{equation}
We let
$$g\lalp\co  (TX)\hhalp_0\llra Q\hhalp_0
$$
be the induced smooth map; we let $\cA_1$ and $\cA_2$ be the sheaves
of smooth functions of $(TX)_0\hhalp$ and $Q\hhalp_0$ respectively.
We claim that there are smooth isomorphisms as shown below that make
the diagram commutative
\begin{equation}\label{20}
\begin{CD} g\lalp\sta\bl\cE\lalp\otimes\cA_2\br @>>>
g\lalp\sta\bl
\Ob\hhalp_Q\otimes \cA_2\br\\
@V{\cong}VV @V{\cong}VV\\
\tilde\cE\lalp\otimes\cA_1 @>>> \Ob\hhalp_{TX}\otimes \cA_1.
\end{CD}
\end{equation}
Indeed, for any open $U_a\sub X$ with trivialization $TX|_{U_a}\cong
V\times U_a$ and open $U_b\sub Gr$ with trivialization
$Q|_{U_b}\cong V\times U_b$ and satisfying $g(U_a)\sub U_b$, the
isomorphism $g_0$ in \eqref{iso4} defines a smooth map $h\mh U_a\to
GL(V)$ that makes the diagram commutative
$$\begin{CD}
Q\times_{Gr}X\ @>{g_0}>> TX\\
@V{\cong}VV @V{\cong}VV\\
V\times U_a @>{(h,1)}>> V\times U_a\\
\end{CD}
$$
The family of automorphisms $h$ then induce smooth diffeomorphism
$$\begin{CD}
V_0\hhalp\times U_a @>{(h\lalp,1)}>> V\hhalp_0\times U_a\\
@V{\cong}VV @V{\cong}VV\\
Q\hhalp_0\times_{Gr}U_a @>{\cong}>> (TX)\hhalp_0\times_X U_a
\end{CD}
$$
and isomorphisms of sheaves
$$h\ualp\co  \cE\lalp|_{Q\hhalp_0\times_{Gr}U_a}\cong
\cF\lalp\otimes\cO_{U_a}\twomapright{} \cF\lalp\otimes\cO_{U_a}\cong
\tilde \cE\lalp|_{(TX)\hhalp_0\times_X U_a},
$$
where the first and the third isomorphisms are induced by the
construction of $\cE\lalp$ and $\tilde\cE\lalp$ in \eqref{quot1}
while the middle one is induced by $h$.

The homomorphisms $h\ualp$ for a covering of $X$ patch together to
form the left smooth isomorphism in \eqref{20}; it makes that
diagram commutative.

This way, by lifting $t\lalp$ to a smooth section in $\cE\lalp$,
then pulling the lifted section back to a smooth section of
$\tilde\cE\lalp$, and finally pushing forward the new section to a
section in the obstruction sheaf $\Ob_{TX}\hhalp$, we obtain a
smooth section. Obviously, the resulting section the original pull
back $\tilde t\lalp$; therefore, $\tilde t\lalp$ is smooth.
Likewise, because $g\mh X\to Gr$ is a smooth embedding, because each
stratum of $C\lalpi$ is submersive onto $Gr$, and because as subsets
$$E\lalp|_{Q\hhalp_0\times_{Gr}X}\supset
D\lalpi\cap \bl Q\hhalp_0\times_{Gr}X\br=\tilde D\lalpi\subset
\tilde E\lalp
$$
under the isomorphism \eqref{20}, we see immediately that a lift of
$\tilde t\lalp$ intersects transversally with the cone cycle $\tilde
C\lalp$ if and only if $Q\hhalp_0\times_{Gr}X$ intersects
transversally with the pseudo-cycle $D(t\lalp)$. But this is
possible if we choose $g$ general. Therefore, we can choose $g$ that
makes all pull back sections $\tilde t\lalp$ satisfy the first
conclusion of \fullref{req}. Since the collection $t\lalp$
satisfies the second conclusion of \fullref{req}, so does the
collection $\tilde t\lalp$.

Apply the same argument to the standard embedding
$$G_N\co  Gr(N,3)\llra Gr(N+1,3)
$$
we see that we can choose sections $t\lalp$'s so that
$$\delta_{[n],Gr(N,3)}=g_N\sta(\delta_{[n],Gr(N+1,3)})
$$
are stable under $G_N$. Therefore, their homology classes
$$[\delta_{[n],Gr(N,3)}]\in H_6(Gr(N,3);\ZZ)
$$
is stable under inclusions.

It remains to use the sections $\tilde t\lalp$ to get smooth
sections $s\lalp$ of $\Ob\hhalp_X$. First, because $\cV\hhalp_0$ is
an open subset of $(TX)\hhalp_0$, the section $\tilde t\lalp$
restricts to a section of the obstruction sheaf $\Ob\hhalp_\cV$ of
$\cV\hhalp_0$. Without much confusion, we shall denote the
restriction section by $\tilde t\lalp$ as well.

Our next step is to take the induced sections $\tilde s\lalp$ of
$\Ob\hhalp_X|_{U\hhalp}$ under the smooth isomorphisms
$$\Psi\lalp\sta\Ob\hhalp_X\cong \Ob\hhalp_\cV
$$
covering the smooth isomorphism $\Phi\lalp\mh \cvalp_0\to
U\hhalp\sub X\hhalp$ in \eqref{14.1} and show that they are smooth
and satisfy the conclusion of \fullref{req}.

We begin with more notations for $X\hhalp$. We pick a locally free
sheaf $\bar\cE\lalp$ on $X\hhalp$ making $\Ob\hhalp_X$ its quotient
sheaf; we let $\bar C\lalp\sub \bar E\lalp$ be their associated
normal cones with irreducible decomposition
$$\bar C\lalp=\sum \bar m\lalpi \bar C\lalpi.
$$
Like before, we denote by $\bar S\lalpi\sub\bar C\lalpi$ its open
stratum and denote by $\bar T\lalpi\sub X\hhalp$ the image of $\bar
S\lalpi$.

Differing from the case of $Q/Gr$, we shall work with the
restriction of the obstruction sheaf to $\bar T\lalpi$. To this end,
we let
$$\bar\cW\lalpi=\Ob\hhalp_X\otimes_{\cO_{X\hhalp}}\cO_{\bar
T\lalpi};
$$
let $\bar W\lalpi$ be the associated vector bundle and let
$$\bar\xi\lalpi\co \bar E\lalp|_{\bar T\lalpi}\llra \bar W\lalpi
$$
be the surjective homomorphism induced by
$\bar\cE\lalp\to\Ob\hhalp_X$.

We shall do the same for $(TX)\hhalp_0$. We let $W\lalpi$ be the
associated vector bundle of the restriction sheaf
$$\tilde\cW\lalpi=\Ob\hhalp_{TX}
\otimes_{\cO_{(TX)\hhalp_0}}\cO_{\tilde T\lalpi},
$$
which is locally free over $\tilde T\lalpi$; we let
$$\tilde\xi\lalpi\co  \tilde E\lalp|_{\tilde T\lalpi}\llra \tilde W\lalpi
$$
be induced by $\tilde\cE\lalp\to\Ob\hhalp_{TX}$.

Our next step is to show that possibly after re-indexing the $i$'s
we have
\begin{equation}\label{21.2}
\Psi\lalp(\bar T\lalpi\cap U\hhalp)=\tilde T\lalpi\cap \cvalp_0;
\end{equation}
and that under the canonical isomorphism $\Psi\lalp\sta \bar
W\lalpi\cong \tilde W\lalpi$ we have
\begin{equation}\label{21.3}
\Psi\lalp\bl \tilde\xi\lalpi(\tilde D\lalpi)\br=\bar\xi\lalpi(\bar
D\lalpi)\and m\lalpi=\bar m\lalpi.
\end{equation}
The proof is straightforward. Let $x\in X$ be any element and let
$\varphi_x\mh\cV_x\to \cU_x\sub X$ be the analytic open embedding
provided by \fullref{lem2.1}. By the construction of the
projection $U\hhalp\to X$, its fiber over $x$ is canonically
isomorphic to
$$U\hhalp\times_{X^n} \cU_x^n\times_{T_xX}0
$$
in which the map $\cU_x^n\to T_xX$ is the composite
$$\varpi_x\co  \cU_x^n\twomapright{\cong}\cV_x^n\twomapright{\sub}
 T_xX^n\twomapright{\text{ave}} T_xX
$$
with the averaging map. Combining the local isomorphism Lemma (\fullref{St2}), the base change property of the obstruction sheaves
\eqref{9.3}, the invariance of the obstruction theory (\fullref{inv5}) and the invariance of normal cone \eqref{cone1}, we see
immediately that to each stratum $S$ of $X\hhalp$, the induced map
$$S\times_{X^n}(\cU_x)^n\twomapright{\varpi_x} T_xX
$$
is a submersion. This shows that the standard stratification of
$X\hhalp$ induces the standard stratification of
$(\cU_x)\hhalp\times_{T_xX}0$. By using $\cU_x\cong \cV_x$ and the
open inclusion $\cV_x\sub T_xX$, the standard stratification of
$(T_xX)\hhalp_0$ also induces the standard stratification of
$(\cU_x)\hhalp\times_{T_xX}0$. Therefore, the restrictions of the
standard stratification of $(T_xX)\hhalp_0$ and of the standard
stratification of $X\hhalp$ to $\cV_0\hhalp\times_X x$ coincide.
Consequently, after re-indexing $\bar C\lalpi$, we will have
\eqref{21.2}.

Then by the base change property of the obstruction sheaves, we
automatically have canonical isomorphism
$$\Psi\lalp\sta \bar W\lalpi\cong \tilde W\lalpi;
$$
applying the invariance results of the obstruction sheaf and of
obstruction theory, we get the identity \eqref{21.3}.

Once we have these, we immediately see that the section $\tilde
s\lalp$ induced by the isomorphism
$$\Psi\lalp\sta(\Ob\hhalp_X)\cong \Ob\hhalp_{TX}|_{\cV_0\hhalp}
$$
is a smooth section of $\Ob\hhalp_X$ over $U\hhalp$; that it
intersects transversally with the normal cone of $\Ob\hhalp_X$ and
the collection $\{\tilde s\lalp\}$ satisfies the conclusions of
\fullref{req} over $U\hhalp$.

To continue, we shall comment on the role of $0<c\ll 1$ and
$0<\eps\ll 1$. In choosing sections $t\lalp$ according to \fullref{req}, we use the smallness of $c$ to force the resulting
pseudo-cycle $\delta_{[n],Q}$ to lie entirely in the $2c$--tubular
neighborhood of the zero section of $Q_0^n/Gr$. As to $\eps$, after
picking $t\lalp$ we use $\eps$ to choose pseudo-cycle
representatives of $D(t\lalp)$ to ensure that the resulting
pseudo-cycle $\delta_{[n],Q}$ represents a homology class in
$H\lsta(Gr,\ZZ)$.

After that, we pick a smooth embedding $g\mh X\to Gr$ to pull the
sections $t\lalp$ back to sections $\tilde t\lalp$ in
$\Ob\hhalp_{TX}$. By choosing $g$ in general position, we can be
sure that the degree of the discrepancy cycle $\delta_{[n],TX}$
constructed using $\tilde t\lalp$ coincides with the pull back of
$\delta_{[n],Q}$ under the induced map $(TX)^n_0\to Q^n_0$. This
time because each $D(\tilde t\lalp)$ is a zero-dimensional
pseudo-cycle, they are cycles automatically. Likewise, because
$\delta_{[n],Q}$ lies in the $2c$--tubular neighborhood of the zero
section of $Q_0^n/Gr$, for $c$ sufficiently small, $\delta_{[n],TX}$
lies entirely in $\cV_0^n$.

The next step is to form the sections $\tilde s\lalp$ of
$\Ob\hhalp_X$ over $U\hhalp$; restrict them to a compact subset
$K\hhalp\sub U\hhalp$ and then extend the restrictions to all
$X\hhalp$ so that the resulting sections $\{s\lalp\}$ satisfy the
conclusion of \fullref{req}. Now let $\Delta_n\sub X^n$ be the top
diagonal $\{(x,\cdots,x)\mid x\in X\}$ as before and let
$\Delta_{n,3c}$ be the $3c$--tubular neighborhood of $\Delta_n$ in
$X^n$. Because $c$ is sufficiently small, we can choose a compact
$K\sub X^n$ so that
$$\Delta_{n,3c}\sub K\sub \cV_0^n\sub X^n.
$$
Thus if we choose $K\hhalp=X\hhalp\times_{X^n}K$, then the
discrepancy cycle $\delta_{[n]}$ constructed using $D(s\lalp)$ is
entirely contained in $\Delta_{n,3c}$; thus is contained in
$U\hhalp$; thus coincide with $\Psi_{[n]}(\delta_{[n],TX})$. Here
the last statement holds because $\delta_{[n],TX}\sub \cV_0^n$.

This completes the proof of \fullref{6.61}.

\section{The case of complex manifolds}

To generalize this theorem to cover all compact, smooth
three-dimensional complex manifolds, we first need to define their
Hilbert schemes of subschemes and their Donalson--Thomas invariants.
The technique developed in this work readily covers the case of
zero-dimensional invariants.

We now construct the Hilbert scheme of points for such complex
manifolds $X$. We shall achieve this goal indirectly by quoting the
invariance results proved in this paper. We first cover $X$ by open
subsets $U\lalp$ such that each is realized as an open subset of
$\At$. By viewing $U\lalp$ as an open subset of $\At$, we define the
Hilbert scheme of $n$--points
$$I_{U\lalp}(0,n)\sub I_{\At}(0,n)
$$
be the open (analytic) subscheme of all $\xi\in I_{\At}(0,n)$ whose
supports lie in $U\lalp$. For any pair $U\lalp$ and $U\lbe$, the
\fullref{St2} ensures that the open subscheme
$$\Hil{U_{\alpha\beta}}\sub \Hil{U\lalp}
$$
is canonically isomorphic to the open subscheme
$$\Hil{U_{\alpha\beta}}\sub \Hil{U\lbe}.
$$
Thus the collection $\Hil{U\lalp}$ glue to form an analytic scheme,
the Hilbert scheme of $n$ points in $X$:
$$\Hil{X}.
$$
The scheme $I_X(0,n)$ comes with the usual obstruction sheaf,
obstruction theory and virtual normal cone. Because of the
invariance results stated or proved in this paper, the obstruction
sheaves $\Ob\lalp$ of $\Hil{U\lalp}$ glue together to form the
obstruction sheaf of $I_X(0,n)$. It can also be defined as the
traceless relative extension sheaf of the universal ideal sheaf of
$\cZ\sub X\times I_X(0,n)$:
$$\Ob=\ext^2_{\pi_2}(\cI_\cZ,\cI_\cZ)_0.
$$
Over each open $\Hil{U\lalp}$, we can find locally free sheaf
$\cE\lalp$ making the obstruction sheaf $\Ob\lalp$ of $\Hil{U\lalp}$
its quotient sheaf; we can also construct its associated normal cone
$C\lalp$ in the associated vector bundle $E\lalp$. Using
$$C\lalp\sub E\lalp\and \cE\lalp\lra\Ob\lalp,
$$
we can define the notion of smooth sections of $\Ob\lalp$ and when
smooth section $s\lalp$ of $\Ob\lalp$ intersects transversally with
the normal come. Because such notion is consistent when restricted
to open subsets $\Hil{U\lalpbe}$, we can make sense of smooth
sections of $\Ob$ on $I_X(0,n)$ and when it intersects transversally
with the normal cone of $\Ob$.

We then define the virtual cycle $I_X(0,n)$ be the homology class
$$[D(s)]\in H_0(I_X(0,n);\ZZ)
$$
represented by the pseudo-cycle $D(s)$ constructed by intersecting
the graph of $s$ with the normal cone in $\Ob$.

To show that the cycle $[D(s)]$ is well-defined, namely it is
independent of the choice of $s$, we need to show that for different
smooth sections $s$ the cycles $D(s)$ are homotopy equivalent. This
is true because we can find a stratification of $I_X(0,n)$ and of
the cone so that each stratum is the complement of finitely many
closed analytic subvarieties in a closed analytic variety.

Combined, this proves

\begin{theo}
Let $X$ be a compact, smooth three dimensional complex manifold.
Then the so constructed Hilbert scheme of points $I_X(0,n)$ has a
well-defined virtual cycle
$$[I_X(0,n)]\virt\in H_0(I_X(0,n);\ZZ)
$$
that is represented by the cycle $D(s)$ after intersecting a smooth
section $s$ transversally with the normal cone in the obstruction
sheaf $\Ob$.
\end{theo}

Once the virtual cycle is constructed, then the proof of this paper
applies line to line to $I_X(0,n)$ to conclude that
$$\deg[I_X(0,n)]\virt
$$
is expressible by the same universal expression in its Chern numbers
as other projective threefolds. Thus

\begin{theo}
The identity
$$\sum_n \deg[I_X(0,n)]\virt q^n=M(-q)^{c_3(T_X\otimes K_X)}
$$
holds for all compact, smooth three dimensional complex manifolds.
\end{theo}

The Hilbert scheme of ideal sheaves of curves for any complex
manifold $X$ can also be defined. In case $X$ has dimension three,
one can also define its virtual cycle and its Donaldson--Thomas
series, along the lines of the work \cite{LT}.

{\bf Acknowledgement}\qua The authors is partially supported by NSF
grants DMS-0200477 and DMS-0244550. 

\bibliographystyle{gtart}
\bibliography{link}

\end{document}